\documentclass[12pt]{iopart}

\usepackage[pagebackref,colorlinks=true,pdfpagemode=none,urlcolor=blue,linkcolor=blue,citecolor=blue]{hyperref}
\usepackage{graphicx}
\usepackage{float}
\usepackage{color}
\usepackage{dsfont}
\usepackage{amssymb, amsthm, cancel, braket}
\usepackage{epstopdf}

\graphicspath{ {/.}{../../Newton-Raphson/}}
\expandafter\let\csname equation*\endcsname\relax
\expandafter\let\csname endequation*\endcsname\relax
\usepackage{amsmath}
\usepackage{enumerate}
\newtheorem{theorem}{\bf Theorem}
\newtheorem{lemma}{\bf Lemma}
\newtheorem{proposition}{\bf Proposition}

\newtheorem{remark}{\bf Remark}
\newtheorem{definition}{\bf Definition}
\newcommand{\R}{\mathds{R}}
\newcommand{\N}{\mathbb{N}}
\newcommand{\A}{\mathcal{A}}
\newcommand{\aaa}{\mbox{$\breve{a}$}}
\newcommand{\fff}{\mbox{$\breve{f}$}}
\newcommand{\dprod}[2]{\left\langle #1, #2 \right\rangle}

\begin{document}

\title[Simultaneous source and attenuation reconstruction in SPECT]{Simultaneous source and attenuation reconstruction in SPECT using ballistic and single scattering data}
\author{M. Courdurier}
\address{Facultad de Matem\'aticas, Pontificia Universidad Cat\'olica de Chile}
\ead{mcourdurier@mat.puc.cl}
\author{F. Monard}
\address{Department of Mathematics, University of Washington, Seattle WA 98195}
\ead{fmonard@math.washington.edu}
\author{A. Osses and F. Romero}
\address{Departamento de Ingenier\'ia Matem\'atica, Universidad de Chile and Centro de Modelamiento Matem\'atico, UMI 2071 CNRS, FCFM, Universidad de Chile}
\ead{fromero@dim.uchile.cl, axosses@dim.uchile.cl}
\bigskip

\begin{abstract} In medical SPECT imaging, we seek to simultaneously obtain the internal radioactive sources and 
the attenuation map using not only ballistic measurements but also first order scattering measurements.
The problem is modeled using the radiative transfer equation by means of an explicit nonlinear operator that gives the ballistic and scattering measurements as a function of the radioactive source and attenuation distributions. First, by differentiating this nonlinear operator we obtain a linearized inverse problem. Then, under regularity hypothesis for the source distribution and attenuation map and considering small attenuations, we rigorously prove 
that the linear operator is invertible and we compute its inverse explicitly. This allows to prove local uniqueness for the
nonlinear inverse problem. Finally, using the previous inversion result for the linear operator, we propose a new type of iterative algorithm for simultaneous source and attenuation recovery for SPECT based on Neumann series and a Newton-Raphson algorithm. 
\end{abstract}


\section{Introduction}
	\subsection{Previous results}

Single-Photon Emission Computed Tomography (SPECT) is a nuclear medicine tomographic imaging technique based on
gamma ray emmision. The idea is to  deliver into a patient a gamma-emitting radioisotope (typically {\em technetium-99m})
that is designed to get attached to certain types of cells or tissues, or distribute in certain region, which then start to emit
gamma rays (for reference see \cite{Atlas}, Chapter 2). This radiation can be measured outside of the patient by a rotating
gamma camera which can identify both the direction and the energy level of the radiation (see Figure \ref{fig1} left). With the information gathered, the goal is to reconstruct the source distribution of the radioisotope inside the patient, hence obtaining
an image of the desired specific tissue in study or the region of interest. SPECT is widely used in monitoring cancer treatment \cite{Keidar2003} and also in neuropsychiatric imaging studies \cite{Quintana2002}.

The mathematical model commonly used to describe the externally measured photons is based on the
Radiative Transfer Equation (RTE) (see e.g. the survey \cite{bal4}) and it requires at least two physical parameters: the radioactive 
source distribution $f$ and the attenuation map $a$. The attenuation map represents the capacity of the medium to 
absorb photons and is, given the medical procedure, an unknown function. As we explained before, the radioactive 
source represents the capacity of the medium to radiate photons, and is the main function to be obtained from SPECT.

The {\em Attenuated Radon transform} (AtRT) plays a central role in SPECT, and particularly in the extensions made in this 
work. The inversion formula for AtRT with known attenuation was obtained independently by Arbuzov, Bukhgeim and
Kazantsev in 1998 \cite{BurkKaz} and by Novikov in 2002 \cite{Novikov} deriving an explicit inverse operator. There are
several generalization for this result, for general geodesics \cite{SaloUhlmann}, for complex valued coefficients \cite{You}
or for more general weight functions \cite{Boman1993, Boman2011,BomStromberg}. There are also invertibility and
stability results for partial measurements. In \cite{Nat1}, injectivity is obtained by measuring in an arbitrarily small open
set of angles. Stability for the direct and inverse problem can be found at \cite{Rullgard180} and inversion of data in 
\cite{Bal180}.

The {\em identification problem} stands in the literature for the simultaneous source and attenuation reconstruction of the pair $(f,a)$ from the AtRT. This problem has been studied for particular cases of attenuation maps: 
constant attenuation (when the 
AtRT reduces to the exponential Radon transform) in \cite{expt2, Solmon1995} (see also \cite{expt1} and the references
therein), radial attenuation in \cite{Puro2013} and piecewise constant attenuation in \cite{Bukhgeim2011}.
For particular cases of source functions, the problem has been also tackled in several papers
\cite{BurkKaz,Bal2,Boman1993,Natterer1981}, and the general non-linear case has been studied in \cite{Stefanov}. Nevertheless, in the general case, examples of non-uniqueness 
then appear: for the weighted Radon transform in \cite{Boman1993}, for the exponential Radon transform in \cite{Solmon1995}, for the non-linear identification problem in \cite{Stefanov} and for the linearized one in \cite{Bal2}. 

Another approach to study the identification problem is by the characterization of the AtRT range. In \cite{Nat2} we can find compatibility properties of the range and in \cite{Novikov2} there is a full characterization of the range. Recently in \cite{Bal2,Stefanov} some local uniqueness and stability results are obtained by using linearization and compatibility conditions for the range. 

The identification problem has also motivated several numerical studies. In many of them \cite{Gourion2002,Ramlau1998,Nat4}, 
the focus is to first obtain a good approximation of the attenuation map instead of treating $(a,f)$ as a pair, called attenuation correction algorithms. For other numerical aspects and reconstructions see for instance \cite{Bronnikov1999, Bronnikov2000,Censor1979,Dicken1999,LuoQianPlamen,Stefanov2014,Manglos1993,Zaidi2003}.

\subsection{Our approach: using lower energy scattering data.}
	   
Our main goal is to reconstruct both the attenuation map and radioactive source distribution of an unknown 
object using the SPECT setting. Although this is the same objective as in the above mentioned identification 
problem, we tackle a different inverse problem by using additional scattering measurements.

\begin{figure}[ht!]
	\centering
	\includegraphics[scale=0.45,trim=0mm 0mm 0mm 0mm,clip]{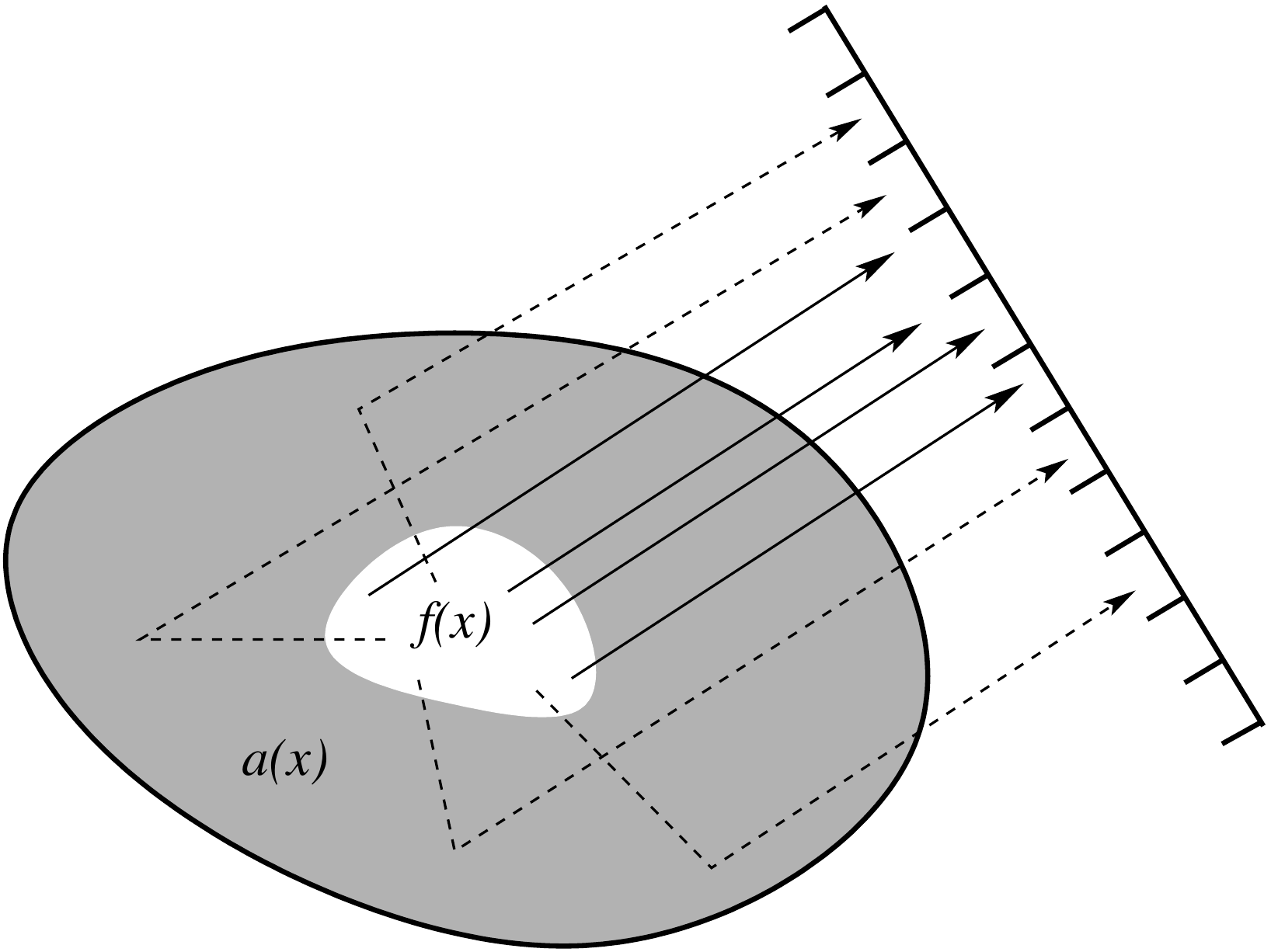}\quad
	\includegraphics[scale=0.55,trim=0mm 0mm 0mm 0mm,clip]{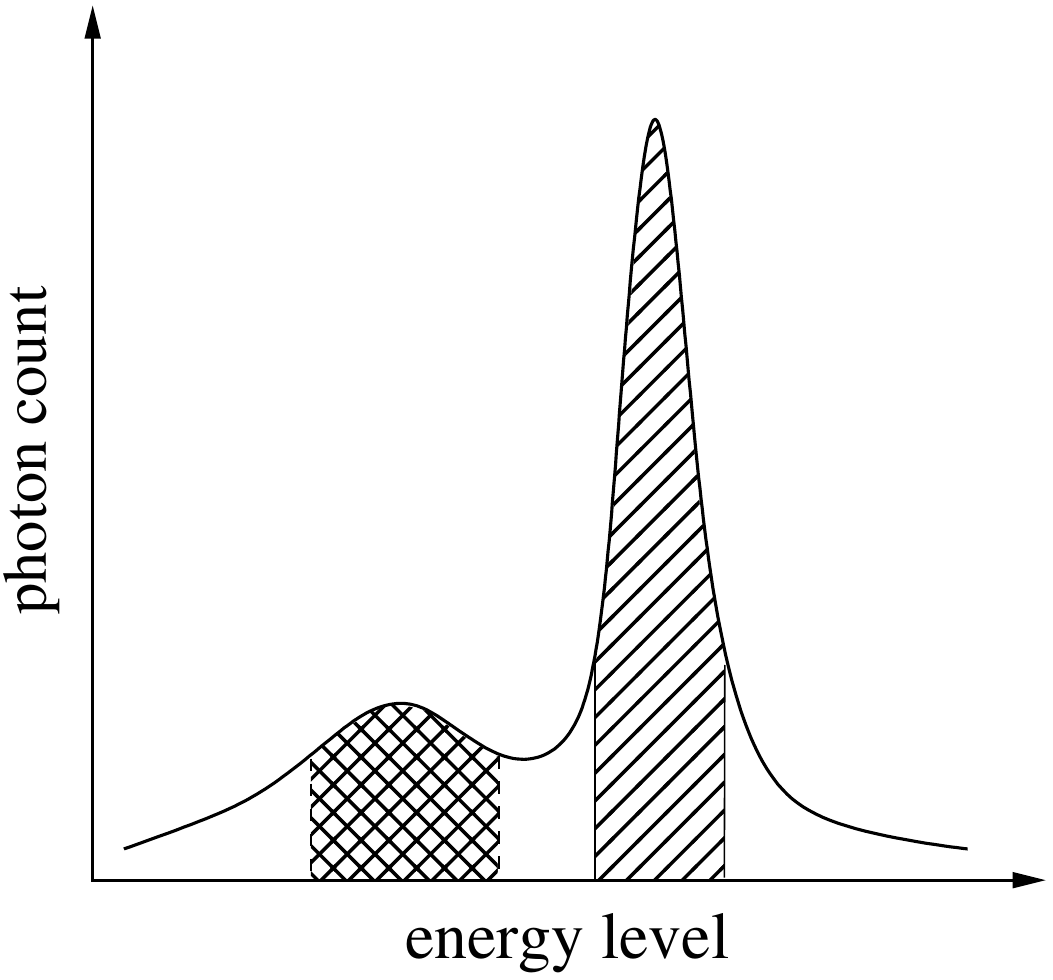}
	\caption{Left: classical SPECT consists in the recovery of the source $f$ by assuming known attenuation $a$ only using high energy ballistic photons (solid arrows). We propose to use available lower energy photon information corresponding mainly to scattered photons (dashed arrows) to recover $f$ and $a$ simultaneously. Right: typical distribution of photons in SPECT showing different energy levels (ballistic: higher energy levels, scattered: lower energy levels).}
	\label{fig1}
\end{figure}

Indeed, we can assume that some additional information can be gathered by measuring scattered 
photons outside the object in study. Each time a photon scatters, it reduces its energy level
(see Figure~\ref{fig1}), and gamma cameras can discriminate the 
energy level of photons. Therefore, we can measure separately the gamma rays exiting the patient that have not scattered (ballistic photons) and the gamma rays exiting the patient that have scattered. Particularly, we are interested 
in measuring photons that just have scattered once (first order scattering photons). 

Considering scattering effects leads to introduce a new unknown coefficient in the model, an scattering coefficient $s(x)$, 
that will describe the scattering behavior of photons inside the object in study. One of our main assumptions will be 
to suppose some relationship between this scattering coefficient and the attenuation coefficient.

In summary, we can assume that we can gather more information using the same standard device and medical 
procedure used for SPECT and without the addition of new technology or other parameters, except for a change
in the protocol for the measurements.

There are three main objectives in this work, the first one is to derive an inverse problem that consider scattering 
effects in the standard mathematical model of SPECT that describes the behavior of photons in a medium, this by 
means of suitable assumptions that allows us to deal with the gathered information from the ballistic and first order scattering photons. The second goal is to reconstruct both the attenuation and source map from the available data. To achieve this goal we study the operator that describes the external measurements by means of a linearization process.
The third objective is to develop an efficient numerical algorithm that use both the ballistic and first order scattering 
photon information to reconstruct both the attenuation and source map of an unknown object.

Notice that the new information given by first order scattering photons traveling in a certain two dimensional 
plane contains information of the whole three dimensional body, because unlike ballistic photons, scattered 
photons do not necessarily travel straight from the source to the gamma camera. Nevertheless, since the 
extension to the three dimensional case is not too different, and in order to simplify notations, we will restrict
the analysis and numerical simulations of this paper to the two dimensional case.

The rest of the paper is organized as follows. In Section 2, we introduce the main notations and definitions and
we develop the mathematical model for the inverse problem. In Section 3, we state the nonlinear and linearized inverse
problems. In Section 4, we first state the main mathematical results of invertibility of the linearized operator (Theorem 
\ref{thrm:yo})
and local uniqueness of the nonlinear inverse problem (Theorem \ref{thrm:nonlinearIP}) and then we present the proofs. Section 5 contains some
numerical experiments that illustrated the feasibility of the proposed SPECT using lower energy scattered photons
in the case of some previously known counterexamples for the identification problem and for other phantoms. 

Natural extensions of this work are: the three dimensional setting, the use of real data, 
the analysis of more general relationships between attenuation and scattering, and the exploration of
alternative numerical reconstructions techniques. We intend to address these points in a forthcoming paper.

\section{Model description}\label{sec:Model}

	\subsection{Notation and functional framework}

Let us introduce the notation of the sets and functional spaces used in this paper. Let $S^1=\{\theta\in \R^2: |\theta|=1\}$ be
the set of directions in $\R^2$, and for $\theta=(\theta_1,\theta_2) \in S^1$, $\theta_1,\theta_2\in\R$,  let
$\theta^\perp=(-\theta_2,\theta_1)$ be its $\pi/2$ counterclockwise rotation.
Let $K$ be a compact set in $\R^2$ of non-empty interior and  let $\tilde{K}$ be a compact set in $\R^2$ slightly larger
than $K$, for simplicity let us consider $K=\{x\in\R^2:|x|\leq 1\}, \tilde{K}=\{x\in\R^2:|x|\leq 2\}$.
For $n\in\N$, $0<\alpha<1$, let $C^\alpha(\R^n)$ be the space of real valued $\alpha$-H\"older continuous functions, for 
$m\in\N\cup\{0,\infty\}$ let $C^m(\R^n)$ be the space of functions with $m$ continuous derivatives,
let $L^2(\R^n)$ be the space of square integrable functions and  let $H^s(\R^n),s>0$ be the classical Sobolev spaces. In
these functional spaces we denote by $||\cdot||_{C^{\alpha}(\R^n)}$, $||\cdot||_{C^{m}(\R^n)}$, $||\cdot||_{L^2(\R^n)}$ 
and $||\cdot||_{H^s(\R^n)}$ the usual norms, omitting $(\R^n)$ in the subscript when the context is clear.
Also let $L^\infty(\R^n)$ be the space of essentially bounded functions, and abusing the notation
let $||\cdot||_\infty$ denote the norm in $L^\infty(\R^n)$ and $C^0(\R^n)$. For $\Omega\subset\R^n$
let $C^{\alpha}(\Omega), C^m(\Omega),  L^2(\Omega), L^\infty(\Omega)$ and $H^s(\Omega)$ be the corresponding 
functional subspaces consisting of functions with support contained in $\Omega$, this may differ from the standard
notation, but it is a convenient notation for this article.

For a function $f:\R^n\times S^1\to\R$ let 
\begin{align*}
&||f||_{C^\alpha(\R^n\times S^1)} = \sup_{\theta\in S^1} ||f(\cdot,\theta)||_{C^\alpha(\R^n)}\\
&\textnormal{and}\\
&||f||_{H^s(\R^n \times S^1)} = \left( \int_{S^1} ||f(\cdot, \phi)||^2_{H^s(\R^n)} d\phi \right)^{1/2},
\end{align*}
and let $C^\alpha(\R^n\times S^1)$ and $H^s(\R^n\times S^1)$ be the smallest Banach spaces
with such norms that contain the compactly supported smooth functions
(these spaces could have more succinctly been described
as $C^\alpha(\R^n\times S^1):= C(S^1; C^\alpha(\R^n))$ and $H^s(\R^n\times S^1):= L^2(S^1; H^s(\R^n))$,
but we adopt the notation commonly used in the related literature, see e.g. \cite{Nat1}).

	\subsection{Integral operators appearing in the model and the inverse problem}

In the modeling and analysis that will be presented in this work there are a number of integral operators that play a crucial role. We proceed to provide a generic definition of these operators,  leaving the discussion of their properties for the next subsection.

For a function $f:\R\to\R$ we let $Hf$ denote its classic Hilbert transform (see e.g. \cite{Duoan2000}). If $g:\R\times S^1\to\R$ then $Hg$ denotes the Hilbert transform of $g(\cdot,\theta)$ for each $\theta\in S^1$. An important integral operator in this paper is the weighted Radon transform.

\begin{definition}[Weighted Radon transform] Let $f:\R^2\to\R$ be a function and $w:\R^2 \times S^1\to\R$ be a weight function, the weighted Radon transform of $f$, with the weight $w$, is defined as,
\begin{align*}
I_w f(s,\theta) = \int_\mathds{R} w(s\theta^\bot+t\theta, \theta) f(s\theta^\bot + t\theta) dt, \quad  s\in\R, \quad \theta\in S^1.
\end{align*}
\end{definition}

We will consider specific weight functions that will themselves be composed of  integral operators, like
 the beam transform.

\begin{definition}[Beam transform] The beam transform of the function $a:\R^2\to\R$, at the point $x \in \R^2$, in the
direction $\theta \in S^1$, is defined as
\begin{align*}
(Ba)(x, \theta) = \int_0^\infty a(x + t\theta) dt,\quad x\in \R^2,\quad \theta \in S^1.
\end{align*}
\end{definition}

The weighted Radon Transform with the exponential of the Beam transform as a weight  is called the attenuated
Radon Transform.

\begin{definition}[Attenuated Radon transform (AtRT)]  \label{atrt}
Let $a, f:\R^2\to\R$, then the attenuated Radon transform of $f$, with attenuation $a$, is defined as
\begin{align*}
R_a f (s,\theta) = \int_\mathds{R}  f(s\theta^\bot+t\theta)e^{-(Ba)(s\theta^\bot+t\theta, \theta)} dt, \quad s \in \R,  \quad \theta \in S^1.
\end{align*}
When $a\equiv 0$ this is called the Radon transform of $f$ and it is denoted as $Rf(s,\theta)$.
\end{definition}

\begin{remark} The attenuated Radon transform is sometimes defined with a different parameterization.
We denote such alternative definition as $R_a^\perp f(s,\theta)$ and it satisfies $R_a^\perp f(s,\theta)=
R_a f(-s,\theta^\perp)$.
This reparameterization does not change the regularity properties of the operator, in particular, results
proved for $R_a$ also apply to $R^\perp_a$.
\end{remark}

For the attenuated Radon transform there exists the following inversion formula (see \cite{Novikov}).

\begin{theorem}[Inverse of the attenuated Radon transform] \label{thrm:iart}
 Let $a,f:\R^2 \to\R$ be continuously differentiable with compact support, then the following complex valued
formula holds pointwise
\begin{align*}
f(x) = \frac{1}{4\pi} \text{\textnormal{ Re}} \text{\textnormal{ div}} \int_{S^1} \theta e^{Ba(x,\theta^\bot)}(e^{-h}He^{h}R_a^\perp f)(x\cdot \theta, \theta)\ d\theta,
\end{align*}
where $h(s,\theta) = \frac{1}{2}(I + iH)R^\perp a(s,\theta)$ (here $i=\sqrt{-1}$, $H$ is the Hilbert transform and $I$ is the
identity operator). For $a$ and $f$ less regular this inversion formula holds in a weaker sense.
\end{theorem}

\begin{definition}\label{def:iart}
 Let $a:\R^2 \to\R$ and $J:\R\times S^1\to \R$, we define $J^\perp(s,\theta)=J(-s,\theta)$ and 
\begin{align*}
R_a^{-1}J (x)  = \frac{1}{4\pi} \text{\textnormal{ Re}} \text{\textnormal{ div}} \int_{S^1} \theta e^{Ba(x,\theta^\bot)}(e^{-h}He^{h}J^\perp)(x\cdot \theta, \theta) d\theta, \quad x\in\R^2,
\end{align*}
where $h(s,\theta) = \frac{1}{2}(I + iH)R^\perp a(s,\theta)$.
\end{definition}

In our analysis we will consider a linearization of the attenuated Radon transform. Such analysis will require to work with a
weight functional of the following form.

\begin{definition} Let $u,v:\R^2\to\R$, we define the weight $w[u,v]:\R^2\times S^1\to\R$ as
\begin{align*}
w[u, v](x,\theta)  &= -\int_{-\infty}^0 e^{-Bu(x+\tau\theta, \theta)}v(x+\tau\theta)\ d \tau, 
 x\in\mathds{R}^2,\quad \theta \in S^1.
\end{align*}
\end{definition}

And the last integral operator that appears in our modeling and analysis is what we call the
focused transform.

\begin{definition}[Focused transform]\label{focused-transform} For $f,a:\R^2\to\R$ we define $M[a,f]:\R^2\to\R$, the
focused transform of the source $f$ with attenuation $a$, as 
\begin{align*}
M[a,f](x) = \int_{S^1} \int_{0}^\infty f(x+t\theta) e^{-\int_0^t a(x+s \theta) ds} dt d\theta,\quad x\in\R^2.
\end{align*}
\end{definition}

	\subsection{Properties of the integral operators and elementary estimates}

The previous integral operators can be defined in different functional spaces with different properties. The 
following result on the continuity of the weighed Radon transform (see e.g. \cite{RullgardEst}) exemplifies the
functional setting in which we will consider the integral operators.

\begin{theorem}\label{thrm:Rullgard} If $1/2 < \alpha \leq 1$, $f\in L^2(K)$ and $w(x,\theta) \in C^\alpha(\R^2 \times S^1)$  then 
\begin{align*}
|| I_w f||_{H^{1/2}(\R \times S^1)} \leq C ||w||_{C^\alpha(\R^2 \times S^1)} || f ||_{L^2(\R^2)},
\end{align*}
where the constant $C$ depends only on the compact set $K$.
\end{theorem}

In order to deal with the different nature of the functional spaces involved, the remainder of this subsection is
devoted to recall some classic results and to provide some technical lemmas that will clarify the computations done
in section \ref{sec:Proofs}. A relationship between Sobolev and H\"older spaces is given by the classic Sobolev embedding.

\begin{theorem}[Sobolev embedding] \label{thrm:sobemb} Let $s,n$ be integers, $\alpha \geq0$, if $(s-\alpha)/n \geq 1/2 $ then 
\begin{align*}
H^s(\R^n) \subset C^{\alpha}(\R^n),
\end{align*}
 and the inclusion is continuous. Also for $s > 1/2$ we have the continuous inclusion
\begin{align*}
H^s(\R) \subset L^{\infty}(\R).
\end{align*}
\end{theorem}

The products of functions that appear will fall under one of the following two lemmas. Their proofs are
obtained by direct calculations.

\begin{lemma} \label{lemma:holdermult}
Let $\Omega \subset \mathds{R}^n$ and $f_1, f_2 \in C^\alpha(\Omega)$, then $f_1 \cdot f_2 \in C^\alpha(\Omega)$ and 
\begin{align*}
||f_1 \cdot f_2||_{C^\alpha(\mathds{R}^n)} & \leq 2 ||f_1||_{C^\alpha(\mathds{R}^n)}||f_2||_{C^\alpha(\mathds{R}^n)}.
\end{align*}
\end{lemma}

\begin{lemma}\label{lemma:hs}
If $f \in H^{1/2}(\mathds{R})$ and $g \in H^s(\mathds{R}),s>1$, then $f\cdot g \in H^{1/2}(\mathds{R})$ and
$$ ||fg||_{H^{1/2}(\mathds{R})} \leq C ||g||_{H^{s}(\mathds{R})} ||f||_{H^{1/2}(\mathds{R})}.$$
\end{lemma}

Next, in Lemmas \ref{lemma:rtrans} and \ref{lemma:btrans}, we recall some of the basic properties of the Radon transform
and the beam transform. These properties follow from their definitions or from results like the projection slice theorem
(see e.g. \cite{Nat1}).

\begin{lemma} \label{lemma:rtrans} We have that
\begin{enumerate}[ a)]
\item If $f(x)=0$ for $|x|>1$ then $Rf(s,\theta)=0, \forall |s|>1, \forall \theta\in S^1$. This is also true for the weighted
Radon transform $I_w f(s,\theta)$.
\item If $f\in C^0(K)$ then $|Rf(s,\theta)|\leq C||f||_\infty, \forall s\in\R,\forall \theta\in S^1$.
\item If $f \in H^t(K), t \geq 0$, then $\forall \theta \in S^1$ the function
$s\mapsto Rf(s,\theta) \in H^t(\R)$ and
\begin{align*}
||Rf(\cdot,\theta)||_{H^t(\R)} \leq C ||f||_{H^t(\R^2)} \  \ \ \forall \theta \in S^1.
\end{align*}
The constants $C$ above only depend on $K$.
\end{enumerate}
\end{lemma}

\begin{lemma} \label{lemma:btrans} Let $a\in C^0(\tilde{K})$ then (recall $\tilde{K}=\{x\in\R^2:|x|\leq2\}$),
\begin{enumerate}[ a)]
\item $\theta \cdot \partial_x B a(x,\theta) = - a (x), \forall x\in\R^2, \forall \theta\in S^1$,
\item $Ba(x,\theta)=0 \text{ if } x\cdot \theta > 2$,
\item $Ba (x,\theta) = R a (x\cdot \theta^\perp,\theta) \text { if } x \cdot \theta <-2$ ,
\item $| Ba(x,\theta) | \leq C ||a||_{\infty}, \forall x\in\R^2, \forall \theta \in S^1$,
\item If $a \in C^\alpha(\tilde{K})$ then $ Ba(x,\theta)\in C^\alpha(\R^2\times S^1)$ and
$||Ba(x,\theta)||_{C^\alpha(R^2\times S^1)}\leq C ||a||_{C^\alpha(\R^2)}.$
\end{enumerate}
The constants above only depend on the compact $\tilde{K}$.
\end{lemma}

We include also the following property for the beam transform.

\begin{lemma} \label{lemma:btrans4}
Let $a\in C^0(\tilde{K})\cap H^1(\tilde{K})$, then
\begin{align*}
|Ba(x,\theta)| \leq C||a||_{H^1(\R^2)}, \forall x\in\R^2, \forall \theta \in S^1.
\end{align*}
\end{lemma}
\begin{proof}
We have
\begin{align*}
|Ba(x,\theta)|  \leq  ||R | a | (\cdot, \theta)||_\infty
 \leq C ||R | a | (\cdot, \theta)||_{H^1(\R)}
 \leq C || ~| a |~||_{H^1(\R^2)}
 \leq C ||a||_{H^1(\R^2)}.
\end{align*}
We used Theorem \ref{thrm:sobemb}, Lemma \ref{lemma:rtrans} and the inequality $|| \, |f| \, ||_{H^1(\R^2)}\leq ||f||_{H^1(\R^2)}$ (e.g. \cite{gilbargtrudinger}).
\end{proof}

And to conclude, we present some technical lemmas on H\"older regularity for functionals that will appear
as or in weight functions.

\begin{lemma} \label{lemma:1}
Let $a \in C^\alpha(\tilde{K})$ and let $k(x,\theta) = e^{-B a (x,\theta^\perp)},  (x,\theta)\in \R^2\times S^1$.
Then $k\in C^\alpha(\R^2\times S^1)$ and 
\begin{align*}
||k||_{C^\alpha(R^2 \times S^1)} \leq C e^{C || a||_\infty }\left(1+||a||_{C^\alpha(\R^2)} \right),
\end{align*}
where $C$ is a constant depending only in the compact set $\tilde{K}$.
\end{lemma}
\begin{proof} Fix $\theta\in S^1$. Since $|k(x,\theta)|\leq e^{||Ba(\cdot,\theta^\perp)||_\infty}$ and
\begin{align*}
|k(x,\theta)-k(y,\theta)|  & = |e^{-B a (x,\theta^\perp)}-e^{-B a (y,\theta^\perp)}|\leq e^{||Ba(\cdot,\theta^\perp)||_\infty}
|Ba(x,\theta^\perp) - Ba(x,\theta^\perp)|,
\end{align*}
we conclude from Lemma \ref{lemma:btrans} that
\begin{align*}
||k(\cdot,\theta)||_{C^\alpha(\R^2)}\leq e^{||Ba(\cdot,\theta^\perp)||_\infty}(1+||Ba(\cdot,\theta^\perp)||_{C^\alpha(\R^2})\leq  C e^{C ||a||_\infty}(1+||a||_{C^\alpha(\R^2}),
\end{align*}
where the constant $C$ is independent of $\theta$ and only depends on $\tilde{K}$.
\end{proof}

We also have the following weight appearing in the inversion of the attenuated Radon Transform.

\begin{lemma} \label{lemma:xeh}
Let $f \in H^2(\tilde{K})$, let $h(s,\theta)=\frac{1}{2}(I + iH)R^\perp f(s,\theta) $ and let  $\varphi\in C^\infty([-2,2])$. Then
$\forall\theta\in S^1$ we have $s\mapsto\varphi(s) e^{h(s, \theta)} \in H^2(\R)$ and
\begin{align*}
|| \varphi(\cdot)e^{\pm h(\cdot,\theta)}||_{H^2(\R)} \leq  C e^{C||f||_\infty } \left(1 + ||f||_{H^2(\R^2)} \right)^2 \ \ \forall \theta \in S^1.
\end{align*}
where $C$ depends only on $\tilde{K}$ and the function $\varphi$.
\end{lemma}
\begin{proof}
From Lemma \ref{lemma:rtrans},
$|| Rf(\cdot, \theta)||_{H^2(\R)} \leq C||f||_{H^2(\R^2)}, \forall \theta\in S^1 $
and since the Hilbert transform is a unitary operator on $H^t(\R), t>0$,
\begin{align*}
|| h(\cdot,\theta)||_{H^2(\R)} \leq C||f||_{H^2(\R^2)}. \forall \theta\in S^1.
\end{align*}
We also have from Lemma \ref{lemma:rtrans},
\begin{align*}
|e^{\pm h(s,\theta)}| = e^{\pm Rf(s,\theta)/2} \leq e^{C ||f||_\infty} \ \ \forall s \in \R,\forall \theta \in S^1.
\end{align*}
With the estimates above, and using the Sobolev inequality $||g||_\infty \leq C||g||_{H^1(\R)}$, it is a direct calculation
to show that the $L^2$ norm of  $\varphi(\cdot)e^{\pm h(\cdot,\theta)}$ and its second derivative, are controlled
as prescribed.
\end{proof}

The focused transform will play an important role in our model and in the inverse problem. The functional framework in which we
will  work with it is the following.

\begin{lemma}  \label{lemma:hold1}
Let $a,f \in C^\alpha(\tilde{K})$, then $M[a,f] \in C^\alpha(\R^2)$ and
\begin{align*}
||M[a, f]||_{C^\alpha} \leq C e^{C ||a||_{\infty}}(1 + ||a||_{C^\alpha}) ||f||_{C^\alpha},
\end{align*}
were $C$ is a constant depending only in the compact set $\tilde K$.
\end{lemma}
\begin{proof} Let us recall that
\begin{align*}
M[a,f](x) = \int_{S^1} \int_{0}^\infty f(x+t\theta) e^{-\int_0^t a(x+s \theta) ds} dt d\theta.
\end{align*}
Let $k_t(x,\theta) = e^{-\int_0^t a(x+s \theta) ds}$, $f_{t\theta}(x)=f(x+t\theta)$ and $\mathds{1}_{\tilde{K}}(x)=1$ if 
$x\in\Tilde{K}$ and $0$ otherwise. Then $||f_{t\theta}||_{C^{\alpha}(\R^2)}=||f||_{C^{\alpha}(\R^2)},\forall t\geq 0,\forall 
\theta\in S^1$, and as in Lemma \ref{lemma:1},
\begin{align*}
||k_t||_{C^\alpha(R^2 \times S^1)} \leq C e^{C || a||_\infty }\left(1+||a||_{C^\alpha(\R^2)} \right),\forall t\geq0.
\end{align*}
Hence
\begin{align*}
||M[a,f]||_\infty &\leq \sup_x \int_{S^1} \int_{0}^\infty ||f_{t\theta}(\cdot) k_t(\cdot,\theta)||_\infty \mathds{1}_{\tilde{K}}(x+t\theta) dt d\theta\\
&\leq  C e^{C || a||_\infty }\left(1+||a||_{C^\alpha(\R^2)} \right)||f||_{C^{\alpha}(\R^2)} \sup_x \int_{S^1} \int_{0}^\infty \mathds{1}_{\tilde{K}}(x+t\theta) dt d\theta\\
&\leq  C e^{C || a||_\infty }\left(1+||a||_{C^\alpha(\R^2)} \right)||f||_{C^{\alpha}(\R^2)}.
\end{align*}
Similarly,
\begin{align*} \left|M[a,f](x)  - M[a,f](y)\right| & \leq \int_{S^1} \int_0^\infty |f_{t\theta}(x) k_t(x,\theta) -
f_{t\theta}(y) k_t(y,\theta)|\\
&\quad\quad\quad\quad\cdot \left( \mathds{1}_{\tilde{K}}(x+t\theta) + \mathds{1}_{\tilde{K}}(y+t\theta) \right) dt d\theta \\
& \leq \int_{S^1} \int_0^\infty 2||f||_{C^\alpha(\R^2)} ||k_t||_{C^\alpha(R^2\times S^1)}|x-y|^\alpha\\
&\quad\quad\quad\quad\cdot  \left( \mathds{1}_{\tilde{K}}(x+t\theta) + \mathds{1}_{\tilde{K}}(y+t\theta) \right) dt d\theta\\
& \leq  C e^{C||a||_\infty} (1 + ||a||_{C^\alpha}) ||f||_{C^\alpha} |x-y|^\alpha,
\end{align*}
where the constant $C$ only depends on the compact set $\tilde{K}$.
\end{proof}
\begin{lemma} \label{lemma:w}
Let $a, f \in C^\alpha(\tilde{K})$ and let $M=M[a,f]\in C^\alpha(\R^2)$. Then $ w[a, f]$ and $ w[a,a \cdot M]\in C^\alpha(\mathds{R} \times S^1)$ and they satisfy
\begin{align*}
||w[a,f]||_{C^\alpha(\R^2 \times S^1)} & \leq C e^{C||a||_{\infty}}(1 + ||a||_{C^\alpha}) ||f||_{C^\alpha}, \\
||w[a, a\cdot M]||_{C^\alpha(R^2 \times S^1)} & \leq C e^{C ||a||_{\infty}}(1 + ||a||_{C^\alpha})^2 
|| a||_{C^\alpha} ||f||_{C^\alpha},
\end{align*}
where $C$ is a constant depending only on the compact $\tilde{K}$.
\end{lemma}
\begin{proof}
We recall that 
\begin{align*}
 w[a,f](x,\theta) = - \int_{-\infty}^0 e^{-\int_0^\infty a(x+t\theta + \tau \theta) d\tau} f(x+t\theta) dt,  x \in \R^2, \theta \in S^1.
\end{align*}
Following the same steps as in the proof of  Lemma \ref{lemma:hold1} we can easily obtain that
\begin{align*}
||w[a, f]||_{C^\alpha(\R^2 \times S^1)} \leq  C e^{C ||a||_{\infty}}(1 + ||a||_{C^\alpha}) ||f||_{C^\alpha}.
\end{align*}
This also implies
\begin{align*}
||w[a, a \cdot M]||_{C^\alpha(\R^2 \times S^1)} \leq  C e^{C ||a||_{\infty}}(1 + ||a||_{C^\alpha})
|| a \cdot  M||_{C^\alpha},
\end{align*}
which together with Lemma \ref{lemma:holdermult} and Lemma \ref{lemma:hold1} concludes the proof.
\end{proof}

	    
	\subsection{Radiative Transfer Equation model and main simplifying hypotheses}
	\label{subsection:simplificatoryHypotesis}

Let $s(x,\theta,\theta^\prime)$ be a scattering kernel that gives us the distribution according to which photons at the spatial point $x \in \mathds{R}^2$, coming from direction $\theta \in S^1$ are scattered in the direction $\theta^\prime \in S^1$. The equation that we use to model the propagation of photons with attenuation $a$, source $f$ and scattering $s$ is, for all $x \in \R^2$ and $\theta \in S^1$
\begin{align}\label{rtesas1}
\begin{aligned} \theta \cdot \nabla_x u(x,\theta) + a(x) u(x,\theta) + \int_{S^1} u(x,\theta)s(x,\theta,\theta^\prime) d\theta^\prime& = f(x) + \int_{S^1} u(x, \theta^{\prime}) s(x, \theta^{\prime}, \theta) d\theta^{\prime} \\
\lim_{t \rightarrow +\infty} u(x - t \theta, \theta)& = 0.
\end{aligned}
\end{align}
The first integral term corresponds to the effect of photons that are scattered away from the path defined by $(x,\theta)$, the second integral term is the opposite, and represents the gamma rays travelling in the spatial point $x\in\R^2$ coming from any direction that by a scattering process take the path defined by $(x,\theta)$.
By introducing the total attenuation:
\begin{equation}
a_T(x) = a(x) + \int_{S^2} s(x,\theta, \theta^{\prime})d\theta^\prime
\end{equation}
then Equation \eqref{rtesas1} can be rewitten as
\begin{equation}
    \theta \cdot \nabla_x u(x,\theta) + a_T(x) u(x,\theta) = f(x) + \int_{S^1} u(x, \theta^{\prime}) s(x, \theta^{\prime}, \theta) d\theta^{\prime}, x \in \R^2, \theta\in S^1. \label{rtesas2}
\end{equation}

Let us introduce $u_i(x,\theta)$ as the intensity of photons that have been scattered $i$ times, thus we can decompose the total intensity $u$ as
\begin{align*}
    u(x,\theta) = \sum_{i=0}^{\infty} u_i(x,\theta),
\end{align*}
(for further reference in this decomposition see e.g. \cite{bal4}), hence Equation \eqref{rtesas2} becomes the
system
\begin{align}\label{rtesas3}
\begin{aligned}
\theta \cdot \nabla_x u_0(x,\theta) + a_T(x)u_0(x,\theta) & = f(x), &\forall x \in \R^2, \theta \in S^1& \\
\theta \cdot \nabla_x u_i(x,\theta) + a_T(x)u_i(x,\theta) & = \int_{S^1} s(x,\theta,\theta^\prime)u_{i-1}(x,\theta^\prime) d\theta^\prime, &\forall i \geq 1, x \in \R^2, \theta \in S^1& \\
\lim_{t \rightarrow +\infty }u_i(x-t\theta,\theta) &=0, &\forall i \geq 0, x \in \R^2, \theta \in S^1. &
\end{aligned}
\end{align}
We first assume isotropy of the scattering kernel $s(x,\theta, \theta^\prime) = s(x, \theta \cdot \theta^\prime)$ i.e. the scattering process just depend on the angle at which photons are scattered, and moreover, we assume we can separate variables for the scattering kernel
\begin{equation}
    s(x,\theta \cdot \theta^\prime) = s(x) k(\theta \cdot \theta^\prime).
\end{equation}
Secondly, we assume that the function $s(x)$ is proportional to the attenuation map (i.e. $\exists C$ such that $Cs(x) = a(x)$) and for the angular variable we assume the scattering kernel is independent of the scattering angle  (i.e. $k(\theta \cdot \theta^\prime)=1/2\pi \quad \forall \theta \cdot \theta^\prime \in [0,1]$). With these assumptions we have for the total attenuation that
\begin{align*}
    a_T(x) = a(x) + \int_{S^2}s(x)k(\theta \cdot \theta^\prime) d \theta^\prime = s(x)(C+1)
\end{align*}
thus redefining $C = (2\pi(1+C))^{-1}$ the system \eqref{rtesas3} becomes 
\begin{align}\label{rtesas4}
\begin{aligned}
\theta \cdot \nabla_x u_0(x,\theta) + a_T(x)u_0(x,\theta) & = f(x),  &\forall x \in \R^2, \theta \in S^1&\\ 
\theta \cdot \nabla_x u_i(x,\theta) + a_T(x)u_i(x,\theta) & = C a_T(x)\int_{S^1}u_{i-1}(x,\theta^\prime) d\theta^\prime, &\forall i \geq 1, x \in \R^2, \theta \in S^1&\\
\lim_{t \rightarrow +\infty } u_i(x-t\theta, \theta) &=0,& \forall i \geq 0,
x \in \R^2, \theta \in S^1.&
\end{aligned}
\end{align}

\begin{proposition} \label{prop:rtesass} If $f$ and $a$ are uniformly line integrable (i.e. $\exists D > 0 : \int_{\R} |f(x+t\theta)| dt \le D, \forall x\in \R^2, \theta \in S^1$) then the system \eqref{rtesas4} has as unique solution  
\begin{align*}
    u_0(x,\theta) &= \int_{-\infty}^0 f(x + t\theta) e^{-\int_t^0 a_T(x+s\theta)ds} dt, \\
    u_i(x,\theta) &= C \int_{-\infty}^0 a_T(x+t \theta)\int_{S^1}u_{i-1}(x+t\theta, \theta^\prime) d\theta^\prime e^{-\int_{t}^0 a_T(x+s\theta)ds}dt, \qquad \forall\ i\ge 1.
\end{align*}
Observe that if $f,a  \in C^0(\tilde{K})$, then they are uniformly line integrable.
\end{proposition}
\begin{proof} The solutions $u_i(x,\theta), i\geq0$ are obtained by direct integration along the characteristics (straight lines) in Equation \eqref{rtesas4}. The line integrability condition ensures by induction that the resulting ODEs can be solved uniquely.			  
\end{proof}

\section{Inverse problem}\label{sec:IP}

\subsection{Measurements and the Inverse Problem}

We assume that the attenuation $a_T(x)$ and the source $f(x)$ are supported in the compact set $K$, representing the
patient. For simplicity from now on we omit the subscript in $a_T$, i.e. the total attenuation is named $a$. For all the
other quantities we keep the notation  of the previous sections.

As measurements we assume that we are able to record $u_0(x,\phi)$, the ballistic photons, and
$u_1(x,\phi)$,
the first order scattering photons, as they exit the patient, i.e. we assume the knowledge of $u_0$ and $u_1$ at all
points outside the support of $a$ and $f$. In summary, the inverse problem that we will study is the reconstruction of the
source and  attenuation maps $f(x)$ and $a(x)$ from the measurement of the ballistic and first order scattering photons
exiting the domain $K$.

Under the hypotheses leading to Proposition \ref{prop:rtesass}, given a source map $f(x)$ and an attenuation coefficient $a(x)$, the intensity of ballistic photons $u_0(x,\theta)$ and the intensity of first order scattering photons 
$u_1(x,\theta)$ at any point $(x,\theta)\in \R^2\times S^1$ is given by 
\begin{align*}
u_0(x,\theta)& = \int_{-\infty}^0 f(x + t\theta) e^{-\int_t^0 a(x+s\theta)ds} dt,  & x\in \mathds{R}^2, \theta \in S^1,& \\
u_1(x,\theta)& = C \int_{-\infty}^0 a(x+t\theta)M[a,f](x+t\theta) e^{-\int_{t}^0 a(x+s\theta)ds}dt,  & x\in \mathds{R}^2, \theta \in S^1,& 
\end{align*}
where $M[a,f](x)=\int_{S^1}u_{0}(x, \theta^\prime) d\theta^\prime$. Therefore, the ballistic and first order scattering photons exiting the domain $K$ correspond to $\A_0, \A_1$, respectively, where we define
\begin{align}
    \mathcal{A}_i (x,\theta)& :=\lim_{\tau\to+\infty}u_i (x+\tau\theta,\theta), \qquad (x,\theta)\in\R\times S^1, \qquad i=0,1.
    \label{eq:A01}
\end{align}
The inverse problem can be rephrased as the reconstruction of $f$ and $a$ on $K$ from knowledge of the Albedo operator 
\[ \mathcal{A}[a,f]=(\mathcal{A}_0,\mathcal{A}_1) = \{(\A_0(x,\theta),\A_1(x,\theta)) , \quad (x,\theta)\in \mathds{R}^2\times S^1\}.   \]
Let us write the operator $\mathcal{A}$ more explicitly. We have
\begin{align*}
\mathcal{A}_0(x,\theta)& = \int_{-\infty}^\infty f(x + t\theta) e^{-\int_t^\infty a(x+s\theta)ds} dt,  &\forall x\in \R^2, \theta \in S^1,& \\
\mathcal{A}_1(x,\theta)& = C \int_{-\infty}^\infty a(x+t\theta)M[a,f](x+t\theta) e^{-\int_{t}^\infty a(x+s\theta)ds}dt,  &\forall x\in \R^2, \theta \in S^1,&\\
M[a,f](x) & = \int_{S^1} \int_{-\infty}^{0} f(x+t\theta) e^{-\int_t^0 a(x+s \theta) ds} dt d\theta,&\forall x\in \R^2, \theta \in S^1,&
\end{align*}
and we observe that $\mathcal{A}_0$ and $\mathcal{A}_1$ are constant along the directed line define by $(x,\theta)$,
i.e. $\mathcal{A}_{i}(x,\theta)=\mathcal{A}_{i}(x+\tau\theta,\theta),\forall \tau\in\R, i=0,1$. Abusing
the notation we write $\forall  s\in \mathds{R}, \theta \in S^1,$
\begin{align*}
\mathcal{A}_0(s,\theta)&:=\mathcal{A}_0(s\theta^\perp,\theta),\\
& = \int_{-\infty}^\infty f(t\theta+ s\theta^\bot) e^{-\int_t^\infty a(\tau\theta+ s\theta^\bot)d\tau} dt,\\
&= R_{a} [f](s,\theta),\\
\mathcal{A}_1(s,\theta)&:=\mathcal{A}_1(s\theta^\perp,\theta),\\
& = C \int_{-\infty}^\infty a( t\theta +s\theta^\bot)M[a,f](t\theta + s\theta^\bot) e^{-\int_{t}^\infty a(\tau\theta+s\theta^\bot)d\tau}dt,\\
&= CR_{a}[ a  M[a,f]](s,\theta).
\end{align*}
Hence we can write the Albedo operator $\mathcal{A}$ as
\begin{align}
\mathcal{A}[a,f]=( R_{a} [f],CR_{a}[ a  M[a,f]]),
\end{align}
and the inverse problem we will study is the inversion of the operator
\begin{align*}
[a,f]\stackrel{\mathcal{A}}{\mapsto}( R_{a} [f],CR_{a}[ a  M[a,f]]).
\end{align*}

\subsection{Formal differential of the Albedo operator and the linearized Inverse Problem}

To study of the invertibility of the Albedo operator $\mathcal{A}$ near
a known source and attenuation pair $(\breve{a},\breve{f})$, supported in $K$, we formally
compute $D\mathcal{A}[\breve{a},\breve{f}](\cdot,\cdot)$ the differential of the Albedo operator
at $(\breve{a},\breve{f})$.

Since $\mathcal{A}[a,f]=( R_{a} [f],CR_{a}[ a  M[a,f]])$ the computation of $D\mathcal{A}$ reduces to the
differentiation of the attenuated Radon transform, which is done in \cite{Stefanov}, and the differentiation of $M[a,f]$.

\begin{proposition}\label{prop:diffA} The formal differential of the Albedo operator at $(\breve{f},\breve{a})$ is
 \begin{align}\label{eqn:lin1}
D\mathcal{A}[\breve{a},\breve{f}](\delta a,\delta f)&= \left( \begin{array}{c}
I_{w[\breve{a}, \breve{f}]}[\delta a] + R_{\breve{a}}[\delta f] \\
 I_{w[\breve{a}, \breve{a} \cdot \breve{M}]} \delta a + R_{\breve{a}}(\delta a \cdot \breve{M})
+ R_{\breve{a}}( \breve{a} \cdot \partial_a \breve{M} \delta a) + R_{\breve{a}}(\breve{a} \cdot M[\breve{a}, \delta f])\end{array}\right).
\end{align}
where 
\begin{align*}
&w[u, v](x,\theta)  = -\int_{-\infty}^0 e^{-Bu(x+\tau\theta, \theta)}v(x+\tau\theta)d \tau, 
 x\in\mathds{R}^2, \theta \in S^1,\\
&\breve{M}=M[\aaa,\fff] \textnormal{ and }\\
&\partial_a \breve{M}\delta a (x) = - \int_{S^1} \int_{0}^\infty \breve{f}( x + t\theta)
 e^{-\int_0^t \breve{a}(x+\tau \theta) d\tau} \int_0^t \delta a(x + s \theta) ds dt d\theta.
\end{align*}
\end{proposition}

\begin{proof}
In \cite{Stefanov} is shown that the formal differential of $(a,f)\mapsto R_{a}[f]$ is $I_{w[a,f]}[\delta a]+R_a[\delta f]$,
which readily implies this result. The computation of $\partial_a\breve{M}\delta a(x)$ is straightforward.
\end{proof}

To study of the operator $D\mathcal{A}[\breve{a},\breve{f}](\delta a,\delta f)$ we
consider some preconditioning. Let us recall that the reference pair $(\breve{a},\breve{f})$
and the perturbation $(\delta a,\delta f)$ are all supported in $K=\{x\in\R^2:|x|\leq1\}$, and also
recall that $I_wf(s,\theta)= R_a f(s,\theta)=0, \forall |s|>1$ if $f$ is supported in $K$.
We will fix $\varphi\in C^\infty([-2,2])$ such that $\varphi(s)=1$ if $|s|\leq 1$ (and $\varphi(s)=0$ if $|s|\geq2$),
and let us define $\chi\in C^\infty(\tilde{K})$ as $\chi(x)=\varphi(|x|)$. The preconditioning of the operator
$D\mathcal{A}[\breve{a},\breve{f}](\delta a,\delta f)$ consist in the following three steps:
\begin{enumerate}[a)]
\item Multiplication component-wise by $\varphi(\cdot)$.
\item Left composition component-wise with $R_{\breve{a}}^{-1}$.
\item Multiplication component-wise by $\chi(\cdot)$.
\end{enumerate}
As we will see later in Proposition \ref{prop:iartreg}, these three steps map continuously
$H^{1/2}(\R\times S^1)$ into $L^2(\R^2)$, depending only on $\breve{a}$. The resulting
preconditioned $D\mathcal{A}[\breve{a},\breve{f}](\delta a,\delta f)$ takes the following
form as a linear operator $(L+Q)$, that we will consider acting on functions supported in $\tilde{K}$.
\begin{definition}\label{def:LQ}
Component-wise we multiply $D\mathcal{A}[\breve{a},\breve{f}](\delta a,\delta f)$ by $\varphi$, we left-compose with 
$R^{-1}_{\breve{a}}$ and then multiply by $\chi$ to obtain the operator $(L+Q)$ defined as
\begin{align*}
L[\breve{a}, \breve{f}](\delta a, \delta f) & =  \left( \begin{array}{c} \chi R_{\breve{a}}^{-1} \varphi I_{w[\breve{a},\breve{f}]}[\delta a] + \delta f \\ \delta a \cdot \breve{M} \end{array} \right), \\
Q[\breve{a}, \breve{f}](\delta a, \delta f) & = \left( \begin{array}{c} 0 \\ \chi R_{\breve{a}}^{-1} \varphi  I_{w[\breve{a}, \breve{a} \cdot \breve{M}]}[ \delta a] + (\aaa \cdot \partial_a \breve{M}\delta a) + (\aaa \cdot M[\aaa, \delta f]) \end{array} \right).
\end{align*}
This operator, originally defined on functions $(\delta a,\delta f)\in L^2(K)\times L^2(K)$, will be extended to be
acting on functions $(\delta a,\delta f)\in L^2(\tilde{K})\times L^2(\tilde{K})$
\end{definition}

The operator $(\delta a, \delta f)\mapsto (L+Q)(\delta a,\delta f)$ is a preconditioned
differential of the Albedo operator $\mathcal{A}$, therefore it represents the linearization of the  originally non-linear
inverse problem when $(\delta a,\delta f)$ are supported in $K$. In section \ref{sec:MainResults} we prove the invertibility of 
the operator $(L+Q)$ in the adequate spaces, providing an explicit inverse and therefore solving the linearized inverse problem.

\subsection{Fr\'echet differentiability of the modified Albedo operator}\label{subsec:Fdiff}

The term $I_{w[\breve{a},\breve{f}]}(\delta a)$ in the differential $D\mathcal{A}[\breve{a},\breve{f}](\delta a,\delta f)$
from the previous subsection, quickly exemplifies why the calculated differential is only formal and not a Fr\'echet differential:
in order to
have the required regularity of $I_{w[\breve{a},\breve{f}]}(\delta a)\in H^{1/2}(\R\times S^1)$ for $\delta a\in L^2(\R^2)$,
we need $w[\breve{a},\breve{f}]\in C^\alpha(\R^2\times S^1)$ for $\alpha>1/2$, which is not going to be the case
for $\breve{a},\breve{f} \in  L^2(\R^2)$. This obstacle can be overcome by considering a modified
Albedo operator for the measurements, one arising from a model in which the attenuation and the source act in a
more regularized way.

Let $\epsilon>0$ and $K_\epsilon=\{x\in\R^2: |x|\leq 1-\epsilon\}$. Let $F_\epsilon:L^2(K_\epsilon)\to H^2(K)$ be an injective
continuous linear operator from $L^2(K_\epsilon)$ into $H^2(K)$ that satisfy $F_\epsilon g \geq 0$ in $K$ if $g\geq 0$ in $K_\epsilon$. For any
$g\in L^2(K_\epsilon)$ write $g_\epsilon:=F_\epsilon(g)\in H^2(K)$ and assume that the transport of photons for a source and
attenuation map $f,a\in L^2(K_\epsilon)$ is instead described by the following modified Radiative Transfer Equation,
\begin{equation}\label{eq:modRTE}
    \phi\cdot\nabla_x u(x,\phi)+a_\epsilon(x) u(x,\phi)=C a_\epsilon(x) \int_{S^1} u(x,\phi')d\phi + f_\epsilon(x), \quad x\in\R^2,\phi\in S^1.
\end{equation}
In this case, the measurements are represented by the following modified Albedo operator.
\begin{definition}\label{def:modAlbedo} We define the modified Albedo operator $\mathcal{A}_\epsilon$ as
\begin{align}
\mathcal{A}_\epsilon[a,f]:=\mathcal{A}[a_\epsilon,f_\epsilon]
=( R_{a_\epsilon} [f_\epsilon],CR_{a_\epsilon}[ a_\epsilon  M[a_\epsilon,f_\epsilon]]).
\end{align}
\end{definition}
On one hand, this modified Albedo operator adds even more assumptions in the model of the measurements.
On the other hand, it has a better behaved functional structure, which can translate into a more robust implementation
in applications. The functional structure of $\mathcal{A}_\epsilon$ is the following.

\begin{proposition}\label{prop:Fdiff}
The operator $\mathcal{A}_\epsilon[a,f]$ is a well defined operator in the spaces
\begin{align*}
\mathcal{A}_\epsilon:L^2(K_\epsilon)\times L^2(K_\epsilon)\to H^{1/2}(R\times S^1)\times H^{1/2}(R\times S^1)
\end{align*}
and is Fr\'echet differentiable at every point $(a,f)\in L^2(K_\epsilon)\times L^2(K_\epsilon)$, with differential 
 \begin{align}
D\mathcal{A}_\epsilon[a,f](\delta a,\delta f)&=
D\mathcal{A}[a_\epsilon,f_\epsilon]((\delta a)_\epsilon,(\delta f)_\epsilon)
\end{align}
where $D\mathcal{A}[\cdot,\cdot](\cdot,\cdot)$ is the formal differential from Proposition \ref{prop:diffA}.
\end{proposition}

\begin{proof}
Since $a,f\in L^2(K_\epsilon)$ then from Sobolev embedding $a_\epsilon, f_\epsilon\in  H^2(K)\subset C^\alpha(K)$ for $\alpha>1/2$. Using Theorem \ref{thrm:Rullgard} and the Lemmas in Section \ref{sec:Model} it follows that
$R_{a_\epsilon} [f_\epsilon]$ and $R_{a_\epsilon}[ a_\epsilon  M[a_\epsilon,f_\epsilon]]$ are in $H^{1/2}(\R\times S^1)$.

For the Fr\'echet differentiability it is enough to show that the reminder term in the first order approximation is quadratic. For both components of the Albedo operator $\mathcal{A}_\epsilon$ we have to study the expansion of the following
generic term
\begin{equation}\label{derivative}
R_{b+\delta b}[g+\delta g]=R_{b}[g]+R_{b}[\delta g]+I_{w[b,g]}[\delta b]+W_1[g]+W_2[\delta g]
\end{equation}
where
$W_1=I_{w_1[b,\delta b]}[g]$ and $W_2=I_{w_2[b,\delta b]}[\delta g]$ with weights
\begin{align*}
w_1[b,\delta b]&= e^{-Bb}(e^{-B\delta b}-1+B\delta b)\\
w_2[b,\delta b]&= e^{-Bb}(e^{-B\delta b}-1).
\end{align*}
For the first component of the Albedo operator we have to take $(b,\delta b,g,\delta g)=(a_\epsilon,(\delta a)_\epsilon,f_\epsilon,(\delta f)_\epsilon)$. From Theorem \ref{thrm:Rullgard}, the Lemmas in Section \ref{sec:Model} and the definition of $F_\epsilon$, we have
\begin{align*}
\|W_1[f_\epsilon]\|_{H^{1/2}(\R\times S^1)}&\le C\|w_1\|_{C^\alpha(\R\times S^1)}\|f_\epsilon\|_{L^2(K)}\\
&\leq Ce^{C\|a_\epsilon\|_\infty}(1+\|a_\epsilon\|_{C^\alpha(K)})\|B(\delta a)_\epsilon\|_{C^0(\R\times S^1)}
\|B(\delta a)_\epsilon\|_{C^\alpha(\R\times S^1)} \|f_\epsilon\|_{L^2(K)}\\
&\leq Ce^{C\|a_\epsilon\|_\infty}(1+\|a_\epsilon\|_{C^\alpha(K)})\|(\delta a)_\epsilon\|_{C^0(K)}
\|(\delta a)_\epsilon\|_{C^\alpha(K)}
\|f_\epsilon\|_{L^2(K)}\\
&\leq Ce^{C\|a_\epsilon\|_{H^2(K)}}(1+\|a_\epsilon\|_{H^2(K)})\|(\delta a)_\epsilon\|^2_{H^2(K)}
\|f_\epsilon\|_{L^2(K)}\\
&\leq Ce^{C\|a\|_{L^2(K_\epsilon)}}(1+\|a\|_{L^2(K_\epsilon)})\|\delta a\|^2_{L^2(K_\epsilon)}
\|f\|_{L^2(K_\epsilon)}\\
&\leq C\|\delta a\|^2_{L^2(K_\epsilon)}
\end{align*}
where the constant $C$ depends on $K, ||a||_{L^2(K_\epsilon)}$ and $||f||_{L^2(K_\epsilon)}$. Similarly
\begin{equation*}
\|W_2[(\delta f)_\epsilon]\|_{H^{1/2}(\R\times S^1)}\le C\|\delta a\|_{L^2(K_\epsilon)}
\|\delta f\|_{L^2(K_\epsilon)}.
\end{equation*}
In other words, for the first component of the Albedo operators, the reminder terms in the first order
approximation are quadratic and therefore it is Fr\'echet differentiable.

For the second component of the Albedo operator let $M_\epsilon=M[a_\epsilon,f_\epsilon]$ and notice that
\begin{equation*}
(a+\delta a)_\epsilon\, M[(a+\delta a)_\epsilon,(f+\delta f)_\epsilon]=h_\epsilon+(\delta h)_\epsilon,
\end{equation*}
where $h_\epsilon$ and $(\delta h)_\epsilon$ are given by
\begin{align*}
h_\epsilon&= a_\epsilon\,M_\epsilon+a_\epsilon\,\partial_aM[a_\epsilon,f_\epsilon](\delta a)_\epsilon
+a_\epsilon\,M[a_\epsilon,(\delta f)_\epsilon]+(\delta a)_\epsilon\,M_\epsilon,\\
(\delta h)_\epsilon&= a_\epsilon\,(M[(a+\delta a)_\epsilon,f_\epsilon]-M_\epsilon
-\partial_aM[a_\epsilon,f_\epsilon](\delta a)_\epsilon)\\ 
&\quad+a_\epsilon\,(M[(a+\delta a)_\epsilon,(\delta f)_\epsilon]-M[a_\epsilon,(\delta f)_\epsilon])\\
&\quad+(\delta a)_\epsilon(M[(a+\delta a)_\epsilon,f_\epsilon]-M_\epsilon)\\
&\quad+(\delta a)_\epsilon\, M[(a+\delta a)_\epsilon,(\delta f)_\epsilon]
\end{align*}
and we expect each term in $(\delta h)_\epsilon$ to be quadratic.

In the generic Equation \eqref{derivative} we have to take $(b,\delta b,g,\delta g)=(a_\epsilon,(\delta a)_\epsilon,h_\epsilon,
(\delta h)_\epsilon)$ and now we have to control the reminder
\begin{equation*}
W=R_{a_\epsilon}[(\delta h)_\epsilon]+I_{w[a_\epsilon,h_\epsilon-a_\epsilon M_\epsilon]}[(\delta a)_\epsilon]
+W_1[h_\epsilon]+W_2[(\delta h)_\epsilon].
\end{equation*}
But since
\begin{align*}
\|h_\epsilon\|_{L^2(K)}&\le C,\\
\|h_\epsilon-a_\epsilon M_\epsilon\|_{L^2(K)}&\le C(\|\delta a\|_{L^2(K_\epsilon)}+\|\delta f\|_{L^2(K_\epsilon)}),\\
\|(\delta h)_\epsilon\|_{L^2(K)}&\le C(\|\delta a\|^2_{L^2(K_\epsilon)}+\|\delta a\|_{L^2(K_\epsilon)}
\|\delta f\|_{L^2(K_\epsilon)}),
\end{align*}
the same argument as before will show that the $H^{1/2}(\R\times S^1)$ norm of the reminder term $W$ satisfy
a quadratic estimate, and therefore the second component of the Albedo operator is also Fr\'echet differentiable.
\end{proof}

\section{Main Results and Proofs}\label{sec:MainResults}
In this section we present the main results and proofs about the linearized inverse problem, establishing the appropriate framework for the problem and concluding with the invertibility of the operator $(L+Q)$ under some assumptions. 
The main idea is straightforward, to prove the invertibility of the linear operator $(L+Q)$ we will show that $L$ is invertible and that $Q$ is a relatively small perturbation. In Section \ref{sec:MainResults}, we present the main results that build up towards the invertibility of $(L+Q)$ and in Section \ref{sec:Proofs}, we present the proof of the main results and the intermediate technical steps.

\subsection{Main Theorems}\label{ssec:MainResults}
First we describe the functional framework in which we study the operator $(L+Q)$.

\begin{proposition}\label{prop:LQinL2} If $\breve{a}\in H^2(\tilde{K})$ and $\breve{f} \in C^\alpha(\tilde{K}), \alpha>1/2,$ then the
operators $L$ and $Q$ from Definition \ref{def:LQ} are well defined in the following spaces
\begin{align*}
L[\breve{a}, \breve{f}], Q[\breve{a}, \breve{f}]&: L^2(\tilde{K}) \times L^2(\tilde{K}) \rightarrow L^2(\tilde{K})\times L^2(\tilde{K}).
\end{align*}
\end{proposition}

The second step is the invertibility of the operator $L$, which is the dominating component of the operator $(L+Q)$.

\begin{proposition} \label{prop:invL}
Let $\breve{a}\in H^2(\tilde{K}), \breve{f} \in C^\alpha(\tilde{K}), \alpha>1/2$. Let $\breve{M}(x)=M[\breve{a},\breve{f}](x)$ and 
assume that $|1/\breve{M}|$ is bounded in $\tilde{K}$. Then the operator $L[\breve{a}, \breve{f}]$ is left-invertible,
with left-inverse
\begin{align}
L^{-1}[\breve{a}, \breve{f}]&:  L^2(\tilde{K})\times L^2(\tilde{K}) \rightarrow L^2(\tilde{K})\times L^2(\tilde{K}),\\
L^{-1}[\breve{a}, \breve{f}] \left( \begin{array}{c} g \\ h \end{array} \right)& = \left( \begin{array}{c} h/\breve{M}\\
 g -  \chi R_{\breve{a}}^{-1} \varphi I_{w[\breve{a},\breve{f}]}[h/\breve{M}] \end{array}\right),
\end{align}
and
\begin{align*}||
L^{-1}[\breve{a},\breve{f}]||\leq  2 +  C(\tilde{K},||\breve{a}||_{H^2})||1/\breve{M}||_{L^\infty(\tilde{K})}(1+||\breve{f}||
_{C^\alpha}),
\end{align*}
where $C(\tilde{K},||\breve{a}||_{H^2})$ is non-decreasing in the norm of $\breve{a}$.
\end{proposition}

The condition $|1/\breve{M}|$ bounded in $\tilde{K}$ is guaranteed to be fulfilled in the following case.

\begin{proposition} \label{prop:posM}
If $\breve{a}\in L^\infty(\tilde{K}),\breve{f}\in C^\alpha(\tilde{K})$, $f\geq0$ and $f\neq 0$, then $1/\breve{M}\in L^\infty(\tilde{K})$ and
\begin{align*}
\breve{M}(x)\geq C e^{-C||\breve{a}||_\infty} \left( \frac{||\breve{f}||_\infty}{|\breve{f}|_{C^\alpha}} \right)^{2/\alpha} ||\breve{f}||_{\infty}, \forall x\in \tilde{K}.
\end{align*}
where $C$ is a constant that only depends on $\textnormal{diam}(\tilde{K})$.
\end{proposition}

The next step is to show that the operator $Q$, i.e. the remainder part of the operator $(L+Q)$, is relatively small
for $\breve{a}$ small. Observe that this is immediate in the critical case $\breve{a}=0$ since then $Q\equiv0$.

\begin{proposition} \label{prop:Qsmall}
Let $\breve{a} \in H^2(\tilde{K})$, $\breve{f} \in C^\alpha(\tilde{K})$ with $\alpha > 1/2$ and $ ||\breve{a}||_{H^2(\mathds{R}^2)} < D$, then 
\begin{align*}
||Q[\breve{a}, \breve{f}]||_{\mathcal{L}(L^2(\tilde{K}),L^2(\tilde{K}))} \leq C(\tilde{K},D)(1+||\breve{f}||_{C^\alpha(\mathds{R}^2)})||\breve{a}||_{H^2(\mathds{R}^2)}.
\end{align*}
 \end{proposition}

Hence, for $\breve{a}$ small, the operator $(L+Q)$ is a small perturbation of an invertible operator, therefore invertible.

\begin{theorem} \label{thrm:yo}
Let  $\breve{a} \in H^2(\tilde{K})$, $\breve{f} \in C^\alpha(\tilde{K})$ with $\alpha > 1/2$, $\breve{f} \geq 0$ and
$\breve{f} \not \equiv 0$. Then $L^{-1}$ and $Q$ are well defined linear operators in the Banach spaces
\begin{align*}
 Q[\breve{a}, \breve{f}], L^{-1}[\breve{a}, \breve{f}]&:  L^2(\tilde{K})\times L^2(\tilde{K}) \rightarrow L^2(\tilde{K})\times L^2(\tilde{K})\
\end{align*}
 and there exists $D > 0$ such that the operator $(L+Q)[\breve{a},\breve{f}]$ defined on $L^2(\tilde{K}) \times L^2(\tilde{K})$ is left invertible for all $\breve{a} \in H^2(K)$ satisfying $||\breve{a}||_{H^2(\mathds{R}^2)} < D$.
The left inverse is given by
\begin{align*}
(L+Q)^{-1}[\breve{a}, \breve{f}] = \sum\limits_{k=0}^{\infty} \left( - (L^{-1}Q)[\breve{a}, \breve{f}] \right)^k \circ L^{-1}[\breve{a},\breve{f}] 
\end{align*}
and
\begin{align*}
||(L+Q)^{-1}[\breve{a}, \breve{f}]|| \leq 2||L^{-1}||.
\end{align*}
The constant $D$ depends only in the compact $\tilde{K}$ and the norms $||f||_\infty, ||f||_{C^\alpha}$.
\end{theorem}

For the non-linear inverse problem with the modified Albedo operator we can provide a local identification result for some
specific perturbations. Using the notation of Subsection \ref{subsec:Fdiff} let $X\subset L^2(K_\epsilon)$ be a set of functions such that if $g\in X$ then
$||g_\epsilon||_{L^2(K)}\geq C ||g||_{L^2(K_\epsilon)}$ for a constant $C>0$ independent of $g\in X$. The
set $X$ can be considered closed under scalarization.

\begin{theorem} \label{thrm:nonlinearIP}
Let $a,a_1,f,f_1\in L^2(K_\epsilon)$ with $f\geq 0$ and $f\not\equiv 0$. Assume that $\delta a= a_1-a\in X$ and
$\delta f=f_1-f\in X$. Then there exist constants $D,\kappa>0$
depending only in the compact set $K$, the operator $F_\epsilon$, the set $X$ and the norms $||f_\epsilon||_\infty, ||f_\epsilon||_{C^\alpha}$ such that if
$||a||_{L^2(K_\epsilon)}\leq D$,  $||\delta a||_{L^2(K_\epsilon)}\leq \kappa$, $||\delta f||_{L^2(K_\epsilon)}\leq \kappa$
and $(a,f)$ produce the same modified measurements as $(a_1,f_a)$, i.e. if 
$\mathcal{A}_\epsilon[a,f]=\mathcal{A}_\epsilon[a_1,f_1]$, then $a=a_1$ and $f=f_1$.
\end{theorem}

\begin{remark}
The requirements on the operator $F_\epsilon$ and the space $X$ arise from technical considerations and it is not
easy to characterize all the pairs $(F_\epsilon,X)$ satisfying the right conditions. Nonetheless, we can provide a
simple example that satisfies all the conditions and that is of interest in applications.

For the operator $F_\epsilon$. Let $\epsilon>0$ be sufficiently small and let $h_\epsilon\in C^\infty(\R^2)$
satisfy $h_\epsilon \geq 0, h_\epsilon\not\equiv 0$ and $h_\epsilon(x)=0$ if $|x|>\epsilon$. Define
$F_\epsilon (g)= h_\epsilon * g$, the convolution between $h_\epsilon$ and $g$, for $g\in L^2(K_\epsilon)$.

For the set $X\subset L^2(K_\epsilon)$. We say that a partition $\mathcal{P}$ of $K_\epsilon$ is in the
family $\mathbf{P}$ if each $P\in \mathcal{P}$ is measurable and contains a ball of radius $2\epsilon$. For a set $A$ let $\chi_A(x)=1$ if $x\in A$ and $\chi_A(x)=0$ otherwise.
Define
\begin{align*}
X=\{g\in L^2(K_\epsilon) : g(x)=\sum_{P\in \mathcal{P}} c_P\chi_P(x), c_P\in \R, \mathcal{P}\in\mathbf{P}\}.
\end{align*}
This pair $(F_\epsilon,X)$ satisfy all the conditions required in the definition of the modified Albedo operator
(Definition \ref{def:modAlbedo}) and in Theorem \ref{thrm:nonlinearIP} above.

\end{remark}

\subsection{Proofs}\label{sec:Proofs}

	This section is devoted to prove the previous results and is organized as follows. We start 
by proving Proposition \ref{prop:posM}, which is a direct computation. The next two steps consist in obtaining estimates
for the operators $I_w$ and $\chi R^{-1}_a \varphi I_w$. We conclude with the analysis of the operators $L$ and $Q$.

\begin{proof}[{\bf Proof of Proposition \ref{prop:posM}}]
Since $\breve{f}\in C^\alpha(\tilde{K})$ and $f\geq 0, f\neq0$, there exists $\overline x \in \tilde{K}$ such that $\breve{f}(\overline x) = ||\breve{f}||_\infty > 0$ . Let
\begin{align*}
A = \{ x \in \tilde{K},\ \breve{f}(x) \geq \frac{||\breve{f}||_\infty}{2} \}.
\end{align*}
We have
\begin{align*} &|\breve{f}(\overline x ) - \breve{f}(y)|  \leq |\breve{f}|_{C^\alpha}|\overline{x}-y|^\alpha \ \ \ \forall y \in A^c \\
\Rightarrow & \frac{||\breve{f}||_\infty}{2}  \leq |\breve{f}|_{C^\alpha} |\overline{x}-y|^\alpha \ \ \ \forall y \in A^c  \\
\Rightarrow & \text{dist}(\overline{x}, A^c) \geq \left( \frac{||\breve{f}||_\infty}{2|\breve{f}|_{C^\alpha}}
 \right)^{1/\alpha}  = R.
\end{align*}
Hence $B(\overline{x}, R) \subset A \subset\tilde{K}$ and $\breve{f}(x) \geq \frac{||\breve{f}||_{\infty} }{2}, \forall x \in B(\overline{x}, R)$.
For $x \in \tilde{K}$
\begin{align*}
M[\breve{a}, \breve{f}](x) & = \int_{S^1} \int_0^\infty f(x + t\theta) e^{-\int_0^t a(x+s\theta)ds} dt d\theta \\
& \geq e^{-\text{diam}(\tilde{K})||\breve{a}||_\infty} \int_{S^1} \int_0^\infty f(x+t\theta) dt d\theta \\
& \geq e^{-\text{diam}(\tilde{K})||\breve{a}||_\infty} \frac{||\breve{f}||_\infty}{2} \int_{S^1} \int_0^\infty 
 \mathds{1}_{B(\overline{x}, R)}(x + t\theta) dt d\theta, \\
& \geq  e^{-\text{diam}(\tilde{K})||\breve{a}||_\infty} \frac{||\breve{f}||_\infty}{2 }\frac{ \pi R^2}{\text{diam}(\tilde{K})}.
\end{align*}
Concluding that 
\begin{align*}
M[\breve{a}, \breve{f}](x)& \geq  \frac{\pi}{2^{1+2/\alpha}\textnormal{diam}(\tilde{K})} e^{-\text{diam}(\tilde{K})||
\breve{a}||_\infty} \left( \frac{||\breve{f}||_\infty}{|\breve{f}|_{C^\alpha}} \right)^{2/\alpha} ||\breve{f}||_{\infty} \ \ \ 
\forall x\in \tilde{K}. 
\end{align*}
\end{proof}

We proceed with some intermediate results needed for the estimates on $I_{w[\breve{a}, \breve{f}]}$ and $I_{w[\breve{a}, \breve{a} \cdot \breve{M}]}$.

\begin{proposition} \label{prop:Iwcomp}
Let $\breve{a},\breve{f} \in C^\alpha(\tilde{K})$ with $\alpha > 1/2$, then  $$I_{w[\breve{a}, \breve{f}]}, I_{w[\breve{a}, \breve{a} \cdot M[\breve{a}, \breve{f}]]}: L^2(\tilde{K}) \rightarrow H^{1/2}(\mathds{R} \times S^1),$$
and for $\delta a \in  L^2(\tilde{K})$,
\begin{align*}
|| I_{w[\breve{a}, \breve{f}]}[\delta a]||_{H^{1/2}(\mathds{R} \times S^1)}
& \leq  C e^{C ||\breve{a}||_{\infty}}\left(1 + ||\breve{a}||_{C^\alpha}\right) ||\breve{f}||_{C^\alpha} ||\delta a ||_{L^2}, \\
|| I_{w[\breve{a}, \breve{a} \cdot M[\breve{a}, \breve{f}]]}[\delta a]||_{H^{1/2}(\mathds{R} \times S^1)}
& \leq C e^{C||\breve{a}||_{\infty}}\left(1 + ||\breve{a}||_{C^\alpha}\right)^2 || \breve{a}||_{C^\alpha}
||\breve{f}||_{C^\alpha}||\delta a||_{L^2}.
\end{align*}
\end{proposition}

\begin{proof}
Follows directly from Theorem \ref{thrm:Rullgard} and Lemma \ref{lemma:w}.
\end{proof}

\begin{proposition} \label{prop:iartreg}
Let $\breve{a} \in H^2(\tilde{K})$ and $J \in H^{1/2}(\mathds{R}\times S^1)$, then 
\begin{align*}
|| \chi R_{\breve{a}}^{-1}[\varphi J ]||_{L^2(\mathds{R}^2)} & \leq C e^{C ||\breve{a}||_\infty}
\left(1+ ||\breve{a}||_{H^2(\mathds{R}^2)} \right)^5   ||J||_{H^{1/2}(\mathds{R} \times S^1)}.
\end{align*}
\end{proposition}
\begin{proof}
Define $\tilde{J}(s,\theta) = \varphi(s) J(-s,\theta^\perp)$. Let $g \in C^\infty(\R^2)$ with compact support and define 
$\tilde{g} = \chi g$. Using  the expression for $R_{\breve{a}}^{-1}$ in Definition \ref{def:iart} we write
\begin{align*}
    \dprod{\chi R_{\breve{a}}^{-1}[ \varphi J ]}{g}_{L^2(\mathds{R}^2)} &
    = \frac{1}{4\pi}\text{Re} \dprod{\text{div} \int_{S^1} \theta e^{(B \breve{a})(x, \theta^\perp)}\left( e^{-h}He^{h} \tilde{J} \right) (x\cdot \theta, \theta) d\theta}{\tilde{g}(x)}_{L^2(\mathds{R}^2)},\\
    &=\frac{1}{4\pi} \text{Re} \dprod{  \int_{S^1} \theta e^{(B \breve{a})(x, \theta^\perp)}\left( e^{-h}He^{h} \tilde{J} \right) (x\cdot \theta, \theta) d\theta}{\nabla \tilde{g}(x) }_{L^2(\mathds{R}^2)}.
\end{align*}
Since $\chi(x)=0$ for $|x|\geq 2 $ (hence $\tilde{g}(x)=0$ for $|x|\geq 2$) and $\varphi(x\cdot \theta/2)=1, \forall \theta\in S^1$ if $|x|\leq 2$, then
\begin{align*}
    \dprod{\chi R_{\breve{a}}^{-1}[ \tilde{J} ]}{ g}_{L^2(\mathds{R}^2)} &= \frac{1}{4\pi} \text{Re} \dprod{ \int_{S^1} \varphi(x\cdot\theta/2) \theta e^{(B \breve{a})(x, \theta^\perp)}\left( e^{-h}He^{h} \tilde{J} \right) (x\cdot \theta, \theta) d\theta}{ \nabla \tilde{g}(x) }_{L^2(\mathds{R}^2)}.
\end{align*}
Defining $F(s,\theta) = \varphi(s/2) e^{-h}He^{h} \tilde{J} (s,\theta)$ we get
\begin{align*}
    \dprod{\chi R_{\breve{a}}^{-1}[ \tilde{J} ]}{g}_{L^2(\mathds{R}^2)} & = 
    \frac{1}{4\pi} \text{Re} \dprod{\int_{S^1} \theta e^{B \breve{a}(x,\theta^\perp)} F(x\cdot \theta, \theta) d\theta}{\nabla \tilde{g}(x)}_{L^2(\R^2)}\\
    &  =  \frac{1}{4\pi} \text{Re} \dprod{ \int_{x\cdot \theta = s} (\theta \cdot \nabla \tilde{g}(x)) e^{(B \breve{a})(x,\theta^\perp)} dl(x)}{F(s,\theta) }_{L^2(\R\times S^1)} \\
    &= \ \frac{1}{4\pi} \text{Re} \dprod{ \int_{x\cdot \theta = s} \theta \cdot \nabla_x \left( \tilde{g}(x) e^{(B \breve{a} )(x,\theta^\perp)} \right) dl(x)}{F(s,\theta)}_{L^2(\mathds{R} \times S^1)}\\
    &\quad - \frac{1}{4\pi}\text{Re}\dprod{ \int_{x \cdot \theta = s} \tilde{g}(x) \theta \cdot \nabla_x e^{(B \breve{a})(x,\theta^\perp)} dl(x)}{F(s,\theta) }_{L^2(\mathds{R} \times S^1)},
\end{align*}
concluding that
\begin{align} 
    \left| \dprod{\chi R_{\breve{a}}^{-1}[ \tilde{J} ]}{ g}_{L^2(\mathds{R}^2)} \right|
    & \leq C \left|  \dprod{ \partial_s \int_{x\cdot \theta = s}  \tilde{g}(x) e^{(B \breve{a} )(x,\theta^\perp)}  dl(x)}{ F(s,\theta) }_{L^2(\mathds{R} \times S^1)}  \right| \label{eq:righthand1} \\
    & \quad + C \left| \dprod{ \int_{x \cdot \theta = s} \tilde{g}(x) \theta \cdot \nabla_x e^{(B \breve{a})(x,\theta^\perp)} dl(x)}{ F(s,\theta) }_{L^2(\mathds{R} \times S^1)} \right|. \label{eq:righthand2}
\end{align}
Let us bound the terms in \eqref{eq:righthand1} and \eqref{eq:righthand2}. For \eqref{eq:righthand1}, let $k(x,\theta^\perp) = e^{(B \breve{a})(x,\theta^\perp)}$, hence
\begin{align*}
    A_1 &= \left| \dprod{ \partial_s \int\limits_{x\cdot \theta = s}\tilde{g}(x) e^{(B \breve{a})(x,\theta^\perp)}dl(x)}{ F(s,\theta) }_{L^2(\R\times S^1)}\right|\\
    &=  \left|\dprod{\partial_s I_{k} \tilde{g}(-s,\theta^\perp)}{ F(s,\theta)}_{L^2(\mathds{R} \times S^1)}\right|\\
&\leq  ||I_{k} \tilde{g} ||_{H^{1/2}(\mathds{R} \times S^1)} ||F ||_{H^{1/2}(\mathds{R} \times S^1)}.
\end{align*}
Using Theorem \ref{thrm:Rullgard} and Lemma \ref{lemma:1} this implies
\begin{align*}
A_1  & \leq  C e^{C ||\breve{a}||_{\infty}} \left( 1 + ||\breve{a}||_{H^2(\mathds{R}^2)} \right) 
||\tilde{g}||_{L^2(\mathds{R}^2)}||F ||_{H^{1/2}(\mathds{R} \times S^1)}.
\end{align*}
To bound  \eqref{eq:righthand2} we observe that from Lemma \ref{lemma:btrans4},
\begin{align*}
|\theta \cdot \nabla_x  e^{B \breve{a} (x,\theta^\perp)} | & = |e^{B \breve{a} (x,\theta^\perp)} B[ \theta \cdot 
\nabla a ](x,\theta) |\\
& \leq C e^{ C || \breve{a} ||_{\infty}} ||\nabla \breve{a} ||_{H^1(\mathds{R}^2)}  \\
& \leq C e^{C ||\breve{a}||_{\infty}} ||\breve{a}||_{H^2(\mathds{R}^2)}, \quad \forall x \in \mathds{R}^2, \theta \in S^1,
\end{align*}
hence
\begin{align*}
    A_2=& \left| \dprod{ \int_{x \cdot \theta = s} \tilde{g}(x) \theta \cdot \nabla_x e^{(B \breve{a})(x,\theta^\perp)} dl(x) }{ F(s,\theta) }_{L^2(\R\times S^1)} \right|\\
    &\leq C e^{C||\breve{a}||_{\infty}} ||\breve{a}||_{H^2(\mathds{R}^2)}
    || R |\tilde{g}| (s,\theta )||_{L^2(\mathds{R}\times S^1)} || F(s,\theta) ||_{L^2(\mathds{R}\times S^1)} \\
    &\leq  C e^{C||\breve{a}||_{\infty}} ||\breve{a}||_{H^2(\mathds{R}^2)} || \tilde{g} ||_{L^2(\mathds{R}^2)} ||F(s,\theta)||_{L^2(\mathds{R} \times S^1)}.
 \end{align*}
In summary, we have obtained the following bound for all $g\in C^\infty(\R^2)$ with compact support, 
\begin{align*}
|\braket{\chi R_{\breve{a}}^{-1} \tilde{J},g}_{L^2(\R^2)} |
& \leq C e^{C ||\breve{a}||_{\infty}}\left( 1+ || \breve{a}||_{H^2} \right)
||F||_{H^{1/2}(\mathds{R} \times S^1)}  ||g||_{L^2}  . \end{align*}
where $F(s,\theta) = \varphi(s/2)e^{-h(s,\theta)}He^{h(s,\theta)} \varphi (s) J(-s,\theta^\perp)$. We complete the proof with 
the following estimate that uses Lemma \ref{lemma:holdermult} and Lemma \ref{lemma:xeh},
\begin{align*}
||F(s,\theta)||^2_{H^{1/2}(\R \times S^1)} & =  \int_{S^1} ||\varphi(\cdot/2) e^{-h(\cdot,\theta)}He^{h(\cdot,\theta)} 
\varphi(\cdot) J(-\cdot, \theta^\perp)||^2_{H^{1/2}(\R)} d\theta  \\
& \leq  \int_{S^1} ||\varphi(\cdot/2) e^{-h(\cdot,\theta)}||^2_{H^2(\R)} || H e^{h(\cdot,\theta)} \varphi(\cdot)
J(-\cdot,\theta^\perp)||^2_{H^{1/2}(\R)} d\theta \\
& \leq  C e^{C ||\breve{a}||_\infty} \left(1 + ||\breve{a}||_{H^2(\R^2)} \right)^4 \int_{S^1}  ||e^{h(\cdot,\theta)} \varphi(s) 
J(-\cdot,\theta^\perp)||^2_{H^{1/2}(\R)} d\theta \\
& \leq   C  e^{C ||\breve{a}||_\infty} \left(1 + ||\breve{a}||_{H^2(\R^2)} \right)^8 \int_{S^1}
||J(-\cdot,\theta^\perp)||^2_{H^{1/2}(\R)} d\theta \\
& =  C  e^{C ||\breve{a}||_\infty} \left(1 + ||\breve{a}||_{H^2(\R^2)} \right)^8 
|| J ||^2_{H^{1/2}(\R \times S^1)}.
\end{align*} 
\end{proof}

\begin{proposition} \label{prop:opbound1}
Let $\breve{a} \in H^2(\tilde{K})$ with $ ||\breve{a}||_{H^2(\R^2)} < D$ and $\breve{f} \in C^\alpha(\tilde{K})$ with
$\alpha > 1/2$, then  
\begin{align*}
\left| \left|\chi R_{\breve{a}}^{-1}\varphi I_{w[\breve{a}, \breve{f}]} \right| \right|_{\mathcal{L}(L^2(\tilde{K}),L^2(\tilde{K}))}
& \leq C(\tilde{K},D) ||\breve{f}||_{C^\alpha(\R^2)} \\
\left| \left|\chi R_{\breve{a}}^{-1} \varphi I_{w[\breve{a}, \breve{a} \cdot M[\breve{a}, \breve{f}]]} \right| \right|_{\mathcal{L}(L^2(\tilde{K}),L^2(\tilde{K}))}
& \leq C(\tilde{K},D)||\breve{a}||_{H^2(\R^2)} ||\breve{f}||_{C^\alpha(\R^2)}, 
\end{align*}
with  $C(\tilde{K},D)$ a constant only depending on $\tilde{K}$ and $D$ and can be taken non-decreasing in $D$.
\end{proposition}
\begin{proof}
This is obtained directly from the estimates in Proposition \ref{prop:Iwcomp} and Proposition \ref{prop:iartreg}.
\end{proof}

\begin{proposition} \label{prop:opbound2}
Let $\breve{a} \in H^2(\tilde{K})$, $\breve{f} \in C^\alpha(\tilde{K})$ with $\alpha > 1/2$ and
$||\breve{a}||_{H^2(\R^2)} < D$, then 
\begin{align*}
Q[\breve{a}, \breve{f}]: L^2(\tilde{K})\times L^2(\tilde{K}) \rightarrow L^2(\tilde{K})\times L^2(\tilde{K}),
\end{align*}
and
\begin{align*}
||Q[\breve{a}, \breve{f}]||_{\mathcal{L}(L^2(\tilde{K})^2,L^2(\tilde{K})^2)} \leq C(\tilde{K},D)
(1+||\breve{f}||_{C^\alpha})||\breve{a}||_{H^2}.
\end{align*}
 \end{proposition}

\begin{proof}
Let $\delta a, \delta f \in L^2(\tilde{K})$, we bound the three terms defining the second component of $Q[\breve{a},\breve{f}](\delta a,\delta f)$. The first term is bounded directly from Proposition \ref{prop:opbound1},
\begin{align*}
|| \chi R_{\breve{a}}^{-1} \varphi I_{w[\breve{a}, \breve{a} \cdot \breve{M}]}[\delta a]||_{L^2(\mathds{R}^2)} \leq C(\tilde{K},D)||\breve{a}||_{H^2(\mathds{R}^2)}||\breve{f}||_{C^\alpha(\mathds{R}^2)}||\delta a||_{L^2(\mathds{R}^2)}.
\end{align*}
The second term $ (\breve{a} \cdot \partial_a \breve{M}\delta a)$ satisfies
\begin{align*}
| \breve{a} \cdot \partial_a \breve{M} \delta a (x)|  = | \breve{a}(x) \int_{S^1} \int_0^\infty \breve{f}(x+t\theta) 
e^{-\int_0^t \breve{a}(x+\tau\theta) d\tau} \int_0^t \delta a(x + s\theta)ds dt d\theta | \\
 \leq  e^{C ||\breve a||_\infty}||\breve{a}||_{\infty}\mathds{1}_{\tilde{K}}(x) \int_{S^1} \int_\R
| \breve{f}(x+ t\theta)|dt \int_{\R} |\delta a (x+ s\theta)| ds  d\theta\\
\leq C  ||\breve{a}||_{\infty} e^{C  ||\breve{a}||_\infty} ||\breve{f}||_{\infty} \int_{S^1} \int_{\R} | \delta a (x+ s\theta)
|  \mathds{1}_{\tilde{K}}(x) \mathds{1}_{\tilde{K}}(x+s\theta) ds d\theta.
\end{align*}
Since $\mathds{1}_{\tilde{K}}(x)\mathds{1}_{\tilde{K}}(x+s\theta)=0$ for $|s|>\textnormal{diam}(\tilde{K}), x\in\R^2$, computing the $L^2$ norm in $x$ gives
\begin{align*}
|| \breve{a} \cdot \partial_a \breve{M} \delta a||_{L^2(\R^2)}  & \leq 
C(\tilde{K}, D) ||\breve{a}||_{\infty} ||\breve{f}||_{\infty} \int_{S^1} \int_{\R} || \delta a ||_{L^2(\R^2)} 
\mathds{1}_{\{|s|\leq \textnormal{diam}(\tilde{K})\}} ds d\theta,\\
&\leq C(\tilde{K}, D) ||\breve{a}||_{\infty} ||\breve{f}||_{\infty} || \delta a ||_{L^2(\R^2)}.
\end{align*}
We can bound the third term $(\breve{a} \cdot M[\breve{a}, \cdot])$ similarly since 
\begin{align*}
| \breve{a} \cdot M[ \breve{a}, \delta f](x) | & = \left| \breve{a}(x) \int_{S^1} \int_0^\infty \delta f(x + t\theta)
e^{-\int_0^\infty \breve{a}(x+ \tau \theta) d\tau} dt d\theta  \right| \\
& \leq e^{C ||\breve{a}||_\infty}||\breve{a}||_\infty \mathds{1}_{\tilde{K}}(x) \int_{S^1} \int_{\R}  |\delta f (x + t \theta)|
\mathds{1}_{\tilde{K}}(x+t\theta )dt d\theta,
\end{align*}
hence
\begin{align*}
|| \breve{a} \cdot M[ \breve{a}, \delta f]||_{L^2(\mathds{R}^2)}& \leq C(\tilde{K},D) ||\breve{a}||_\infty ||\delta f ||_{L^2(\mathds{R}^2)}.
\end{align*}
These three estimates readily imply the result.
\end{proof}

We have all the estimates needed to prove the main results estated in the previous subsection.

\begin{proof}[{\bf Proof of Proposition \ref{prop:LQinL2}}]
The fact that $Q[\breve{a},\breve{f}]:L^2(\tilde{K})\times L^2(\tilde{K})\to
 L^2(\tilde{K})\times L^2(\tilde{K})$ is established in Proposition \ref{prop:opbound2}. The fact that
$L[\breve{a},\breve{f}]:L^2(\tilde{K})\times L^2(\tilde{K})\to 
L^2(\tilde{K})\times L^2(\tilde{K})$ is a direct consequence of Proposition \ref{prop:opbound1} and the fact that
$\breve{M}=M[\breve{a},\breve{f}]\in L^\infty(\tilde{K})$ by Lemma \ref{lemma:hold1}.
\end{proof}

\begin{proof}[{\bf Proof of Proposition \ref{prop:Qsmall}}]
This is established in Proposition \ref{prop:opbound2}.
\end{proof}

\begin{proof}[{\bf Proof of Proposition \ref{prop:invL}}]
 Given $g,h \in L^2(\tilde{K})$, by the definition of $L^{-1}$,  
\begin{align*}
\left|\left|L^{-1}[\breve{a}, \breve{f}]\left( \begin{array}{c} g \\ h \end{array} \right)\right|\right|^2_{L^2\times L^2}  \leq 
|| h/\breve{M}||^2_{L^2} + 2||g||^2_{L^2} + 2||\chi R_{\breve{a}}^{-1} \varphi I_{w[\breve{a},
\breve{f}]}[h/\breve{M}]||^2_{L^2}.
\end{align*}
By Proposition \ref{prop:opbound1}, for $D=||\breve{a}||_{H^2}$, 
\begin{align*}
||\chi R_{\breve{a}}^{-1} \varphi I_{w[\breve{a}, \breve{f}]}[h/\breve{M}] ||_{L^2} \leq C(\tilde{K},D)
||\breve{f}||_{C^\alpha} ||h/\breve{M}||_{L^2}.
\end{align*}
hence 
\begin{align*}
\left|\left|L^{-1}[\breve{a}, \breve{f}]\left( \begin{array}{c} g \\ h \end{array} \right)\right|\right|^2_{L^2\times L^2}  \leq 
2||g||^2_{L^2} + C(\tilde{K},D)(1+||\breve{f}||_{C^\alpha})^2 ||1/\breve{M}||^2_{L^\infty(\tilde{K})} || h||^2_{L^2},
\end{align*}
proving the proposition.
 \end{proof}

 \begin{proof}[{\bf Proof of Theorem \ref{thrm:yo}}]
From Proposition \ref{prop:invL}, Proposition \ref{prop:posM} and Proposition
\ref{prop:Qsmall},
\begin{align*}
||L^{-1}Q||&\leq||L^{-1}||\cdot ||Q||\\
&\leq \left(2 +  C(\tilde{K},D)(1 + ||\breve{f}||_{C^\alpha}) \frac{||\breve{f}||^{2/\alpha}_{C^\alpha}}{||\breve{f}||^{1+2/\alpha}_{\infty}}\right)\left( C(\tilde{K},D)(1+||\breve{f}||_{C^\alpha})||\breve{a}||_{H^2}\right)\\
&\leq D \left( C(\tilde{K},D)(1+||\breve{f}||_{C^\alpha})^2
(1+\frac{||\breve{f}||^{2/\alpha}_{C^\alpha}}{||\breve{f}||^{1+2/\alpha}_{\infty}})\right)\\
&\leq\frac{1}{2}
\end{align*}
for $D$ sufficiently small, since $C(\tilde{K},D)\geq 0$ is non-decreasing in $D$. Hence $(I+L^{-1}Q)$ is invertible and its inverse can be written as a Neumann series. Therefore $(L+Q)$ is left invertible and its left inverse can be written as
\begin{align*}
(L+Q)^{-1}&=(I+L^{-1}Q)^{-1}L^{-1}\\
&=\sum\limits_{k=0}^{\infty} \left( - (L^{-1}Q) \right)^k L^{-1} .
\end{align*}
From the estimate for $||L^{-1}Q||$ we also get $||(L+Q)^{-1}||\leq 2||L^{-1}||$.
\end{proof}

\begin{proof}[{\bf Proof of Theorem \ref{thrm:nonlinearIP}}]
Assume that $(a,f), (a_1,f_1) \in L^2(K_\epsilon)\times L^2(K_\epsilon)$ satisfy the hypothesis of the theorem.
In particular assume that $\mathcal{A}_\epsilon(a,f)=\mathcal{A}_\epsilon(a_1,f_1)$.
Write $\delta a=a_1-a$ and $\delta f = f_1-f$, from Proposition \ref{prop:Fdiff}  we have
\begin{align*}
\mathcal{A}_\epsilon(a+\delta a,f+\delta f)=\mathcal{A}_\epsilon(a,f)+
D\mathcal{A}[a_\epsilon,f_\epsilon]((\delta a)_\epsilon,(\delta f)_\epsilon)-R
\end{align*}
where $R$ satisfies
\begin{align*}
||R||_{H^{1/2}(\R\times S^1)} \leq C (||\delta a||_{ L^2(K_\epsilon)}+||\delta f||_{ L^2(K_\epsilon)})^2.
\end{align*}
Since the value of the Albedo operator at $(a,f)$ and $(a_1,f_1)$ agree, we have
\begin{align*}
D\mathcal{A}[a_\epsilon,f_\epsilon]((\delta a)_\epsilon,(\delta f)_\epsilon)=R.
\end{align*}
In the relationship above we apply the preconditioning steps before Definition \ref{def:LQ}: component-wise multiply by
$\varphi(s)$, compute the inverse of the attenuated Radon transform with attenuation $a_\epsilon$,
then multiply by $\chi(x)$. We obtain
\begin{align*}
(L+Q)[a_\epsilon,f_\epsilon]((\delta a)_\epsilon,(\delta f)_\epsilon)= \chi R^{-1}_{a_\epsilon} [\varphi R],
\end{align*}
hence
\begin{align*}
((\delta a)_\epsilon,(\delta f)_\epsilon)= (L+Q)^{-1}[a_\epsilon,f_\epsilon](\chi R^{-1}_{a_\epsilon} [\varphi R]).
\end{align*}
The hypotheses for Theorem \ref{thrm:yo} and Proposition \ref{prop:iartreg} are satisfied, therefore
\begin{align*}
||((\delta a)_\epsilon,(\delta f)_\epsilon)||_{L^2(K)\times L^2(K)}
\leq C (||\delta a||_{L^2(K_\epsilon)}+||\delta f||_{L^2(K_\epsilon)})^2.
\end{align*}
Since $\delta a,\delta f\in X$ and $||\delta a||_{L^2(K)}\leq \kappa, ||\delta f||_{L^2(K)}\leq \kappa$ then 
\begin{align*}
||\delta a||_{L^2(K_\epsilon)}+||\delta f||_{L^2(K_\epsilon)}
\leq C\kappa (||\delta a||_{L^2(K_\epsilon)}+||\delta f||_{L^2(K_\epsilon)}).
\end{align*}
For $\kappa$ small enough this implies $\delta a=0$ and $\delta f=0$.
\end{proof}

\section{Numerical experiments}

We now present a {\tt MatLab} implementation of the Newton-Raphson algorithm based on the linearized inverse problem. 

The computational domain is the unit square $[-1,1]^2$ discretized into an equispaced cartesian grid of size $N\times N$ with $N=256$. The quantities of interest $(a,f)$ are supported inside the unit disc $D = \{x^2 + y^2 <1\}$. The computation of the forward measurement operator consists in computing the ballistic and single scattering parts (resp. $\A_0[a,f]$ and $\A_1[a,f]$ as defined in \eqref{eq:A01}, call these measurements $\A$), outgoing traces of the solutions $u_0, u_1$ of system \eqref{rtesas3}.
Such a method is referred to as the {\em iterated source} method (see e.g. \cite{Bal2010}) for solving \eqref{rtesas2}, in the exact same way that in system \eqref{rtesas3}, the term $u_{i-1}$ yields a source term for a transport equation satisfied by $u_i$. Computing such quantities is based on discretizing $S^1$ uniformly and, for each $\theta$ in this discretization, integrating first-order ODEs along lines of fixed direction $\theta$. The latter task is done by computing rotated versions of the map one desires to integrate (e.g. $a$ or $f$) so that the direction of integration coincides with one of the cartesian axes of the image, and the integration along each row is done via cumulated sums, see \cite{Bal2010} for details. In the present case, computing the values of a rotated image is achieved via bilinear interpolations.

Subsequently, the iterative inversion is done by implementing the modified Newton-Raphson scheme 
\begin{align} 
    \begin{split}
	(a^0, f^0) &= (0, 1), \\
	(a^{n+1}, f^{n+1}) &= (a^n, f^n) -  L^{-1} \left( \sum_{k=0}^\infty (-QL^{-1})^k[F_\epsilon(a^n,f^n)]\right) \chi R^{-1}_{a_\epsilon} ( \mathcal{A}_\epsilon(a^n, f^n) - \A),	
    \end{split}
    \label{eq:NRscheme}
\end{align}
where the operators $L,Q$ come from Definition \ref{def:LQ} and the modified albedo operators $\A_\epsilon$ comes from Definition \ref{def:modAlbedo}. 

\begin{remark}
    This algorithm is a modified Newton-Raphson algorithm in the sense that we use the inverse of $D\mathcal{A} [F_\epsilon(a^n,f^n)]$ instead of the inverse of $D\mathcal{A}_\epsilon[F_\epsilon(a^n,f^n)]$ in the right hand side of equation \eqref{eq:NRscheme}. The operator $F_\epsilon$ is introduced in Section \ref{sec:IP} for theoretical purposes and it can be chosen to be an approximation of identity (mollifier). Numerically, doing so adds robustness to the scheme and does not affect the convergence of the algorithm to the correct target functions.
\end{remark}

In all experiments below, 8 iterations of the scheme \eqref{eq:NRscheme} are enough to ensure convergence, and the Neumann series $\sum_{k=0}^\infty (-QL^{-1})^k$ is approximated by its first 4 terms. The implementation of $L^{-1}$ and $Q$ is straightfoward via rotations, cumulated sums and pointwise multiplications/division on the cartesian grid. 

{\bf Axes on figures.} In Figures \ref{fig:pairs} through \ref{fig:disconterrors}, functions of $(x,y)$ are represented on the unit square $[-1,1]^2$. For $i=0,1$, the measurement data $\A_i[a,f]$ (e.g. on Fig. \ref{fig:pairs}, bottom row) are represented by their values for $\A_i(s\theta^\perp,\theta)$ for $\theta\in [0,2\pi]$ on the horizontal axis and $s\in [-1,1]$ on the vertical axis.

In the sections below, we present two series of experiments. Section \ref{sec:nonunique} aims at showing that considering the data $\A_1(a,f)$ in addition to $\A_0(a,f)$ tremendously improves the conditioning of an inverse problem referred to as the identification problem. Secondly, while Section \ref{sec:nonunique} treats the reconstruction of smooth pairs $(a,f)$, Section \ref{sec:nonsmooth} illustrates the performance of our algorithm in the case of discontinuous unknown coefficients from measurements with different levels of noise. 

\subsection{Non-unique pairs and trapping geometries} \label{sec:nonunique}
\begin{figure}[htpb]
  \centering
  \includegraphics[trim = 20 0 20 0, clip, width=0.23\textwidth]{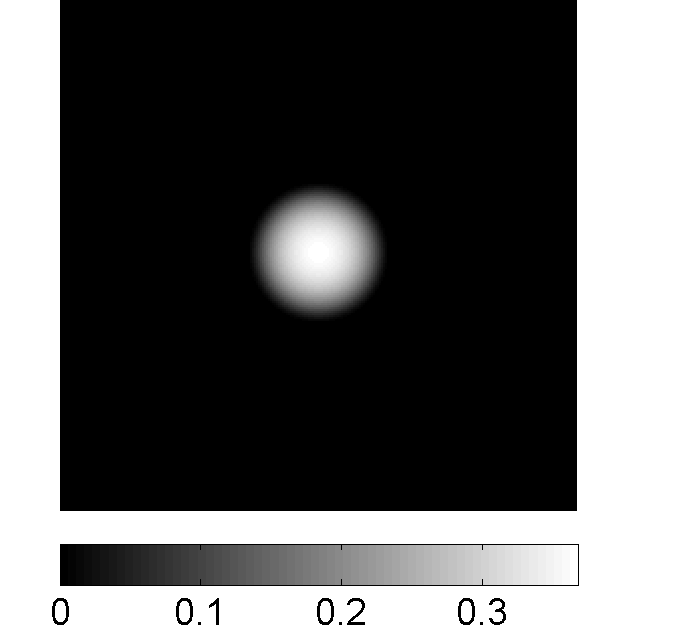}
  \includegraphics[trim = 20 0 20 0, clip, width=0.23\textwidth]{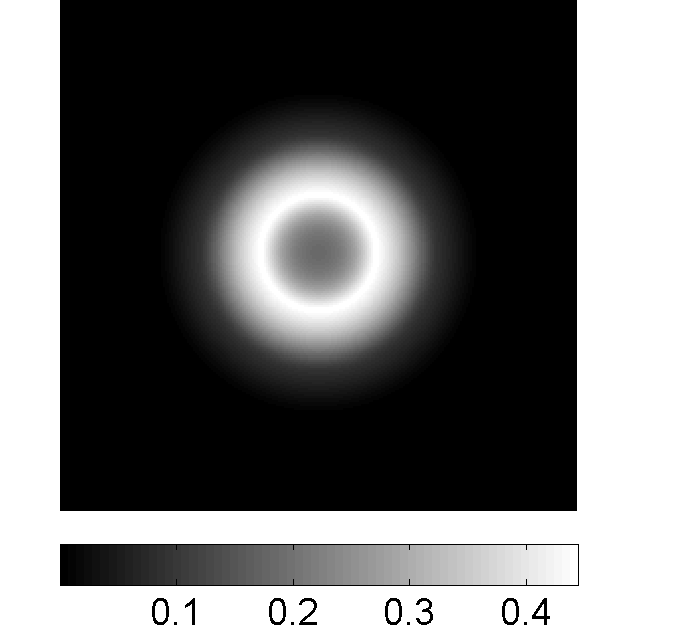}\;\;\;\; 
  \includegraphics[trim = 20 0 20 0, clip, width=0.23\textwidth]{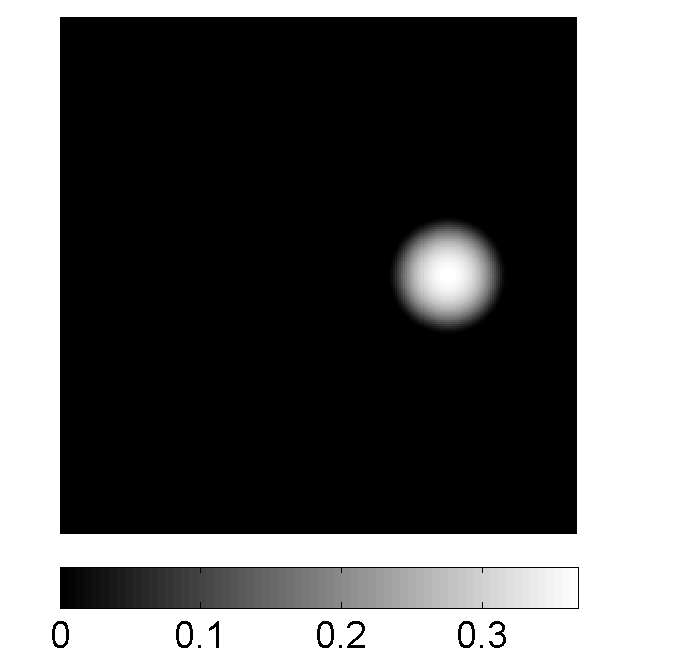}
  \includegraphics[trim = 20 0 20 0, clip, width=0.23\textwidth]{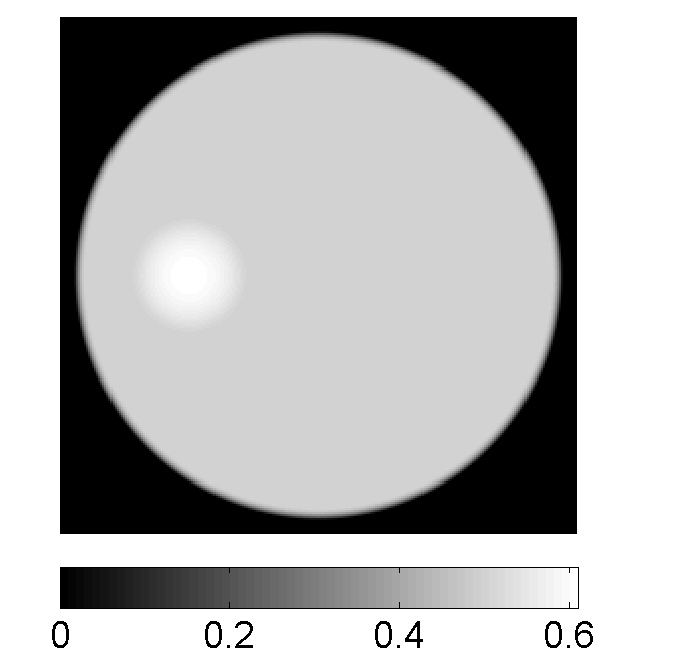} \\
  \includegraphics[trim = 10 0 10 0, clip, width=0.23\textwidth]{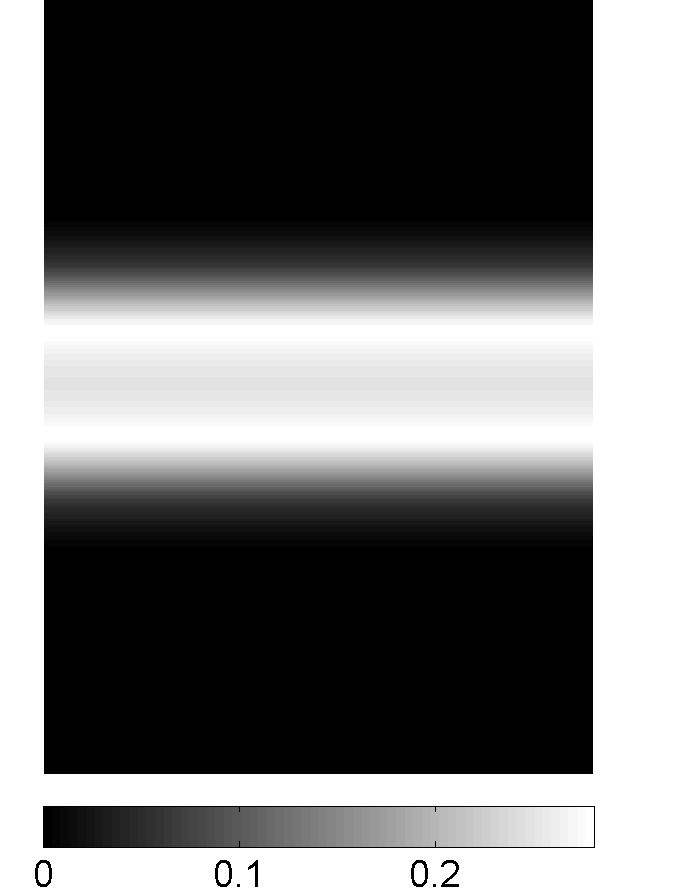}
  \includegraphics[trim = 10 0 10 0, width=0.23\textwidth]{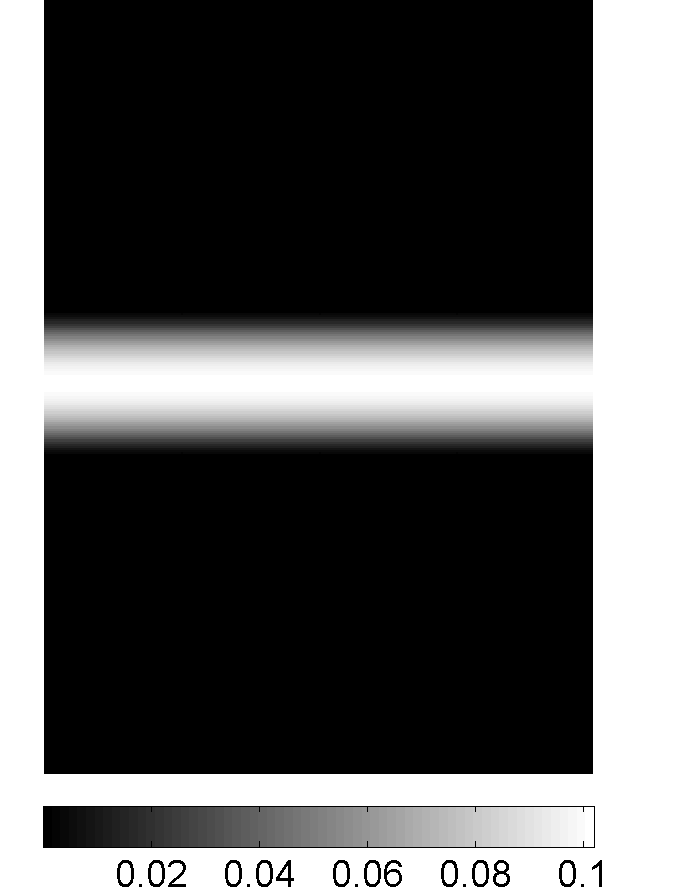}\;\;\;\;
  \includegraphics[trim = 10 0 10 0, width=0.23\textwidth]{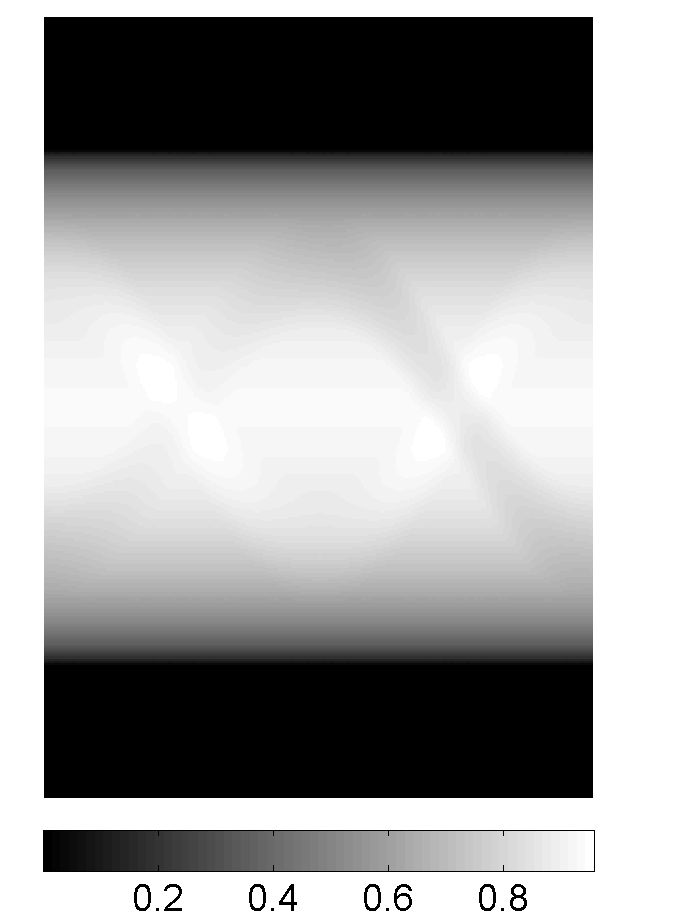}
  \includegraphics[trim = 10 0 10 0, width=0.23\textwidth]{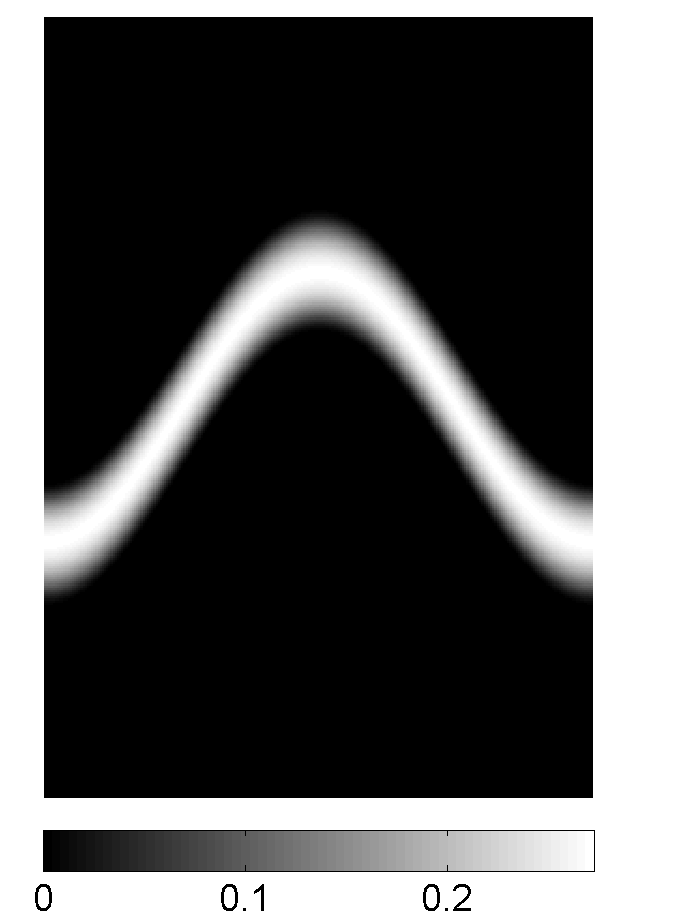}
  \caption{Top row: a ``non-unique'' radial pair $(a_1,f_1)$ (left) and a non-stable pair $(a_2,f_2)$ as defined and studied in \cite{LuoQianPlamen} ($(a_2,f_2)$ are rotated by 90 degrees). Bottom row: the data $(\A_0(a_1,f_1), \A_1(a_1,f_1))$ (left) and $(\A_0(a_2,f_2), \A_1(a_2,f_2))$ (right). }
  \label{fig:pairs}
\end{figure}

The problem of reconstructing the pair $(a,f)$ from only the measurements $\A_0(a,f)$ is known as the {\em identification problem}. Recent theoretical work on the identification problem \cite{Stefanov} and corresponding numerical experiments in \cite{LuoQianPlamen} show that this problem is badly conditioned in at least two ways: 
\begin{description}
  \item[Lack of injectivity.] If $a$ and $f$ are both radial (and smooth enough), there exists a radial function $f_0$ such that $\A_0(a,f) = \A_0 (0,f_0)$. This lack of injectivity prevents some experiments done in \cite{LuoQianPlamen} from converging to the right unknowns.
  \item[Instability.] On the linearized problem (say, reconstruct $(\delta a,\delta f)$ from $\delta \A_0$ around a background $(\aaa,\fff)$), microlocal stability is lost when a certain Hamiltonian flow related to the background $(\aaa,\fff)$ has trapped integral curves inside the domain of interest (referred to as a {\em trapping} geometry). In this case, experiments done in \cite{LuoQianPlamen} show the presence of artifacts in reconstructions. 
\end{description}
The numerical experiments of this section aim at showing that accessing the additional measurement $\A_1(a,f)$ helps at successfully reconstructing both unknowns $(a,f)$ in both scenarios described above. Figure \ref{fig:pairs} displays a pair $(a,f)$ corresponding to each scenario, as well as the corresponding forward data.

\begin{figure}[htpb]
  \centering
  \includegraphics[trim = 20 0 20 0, clip, width=0.23\textwidth]{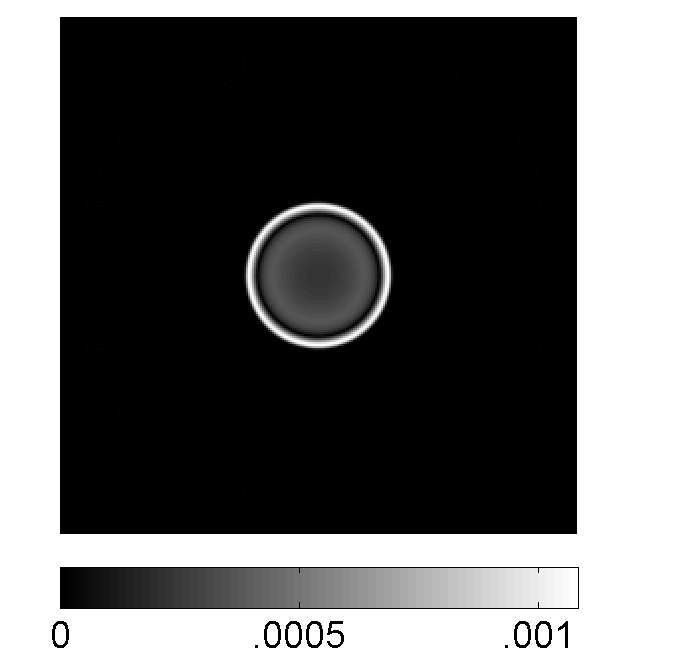}
  \includegraphics[trim = 20 0 20 0, clip, width=0.23\textwidth]{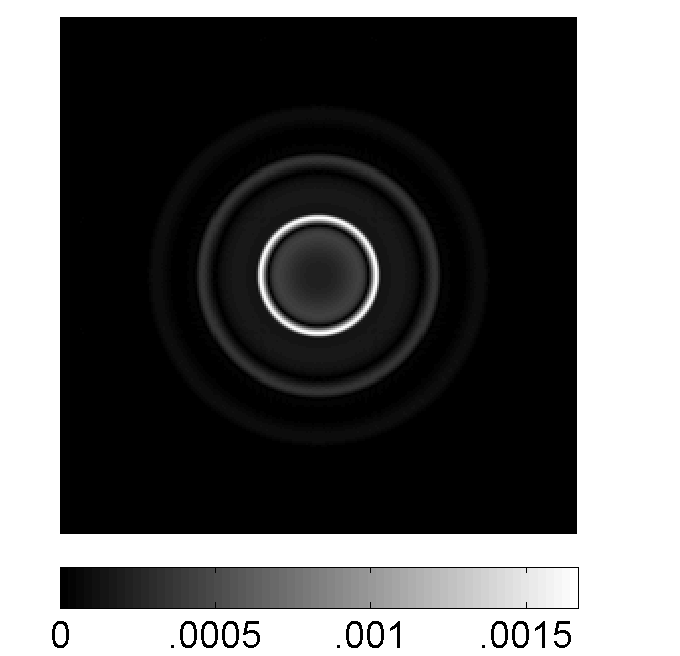}\;\;\;\;
  \includegraphics[trim = 20 0 20 0, clip, width=0.23\textwidth]{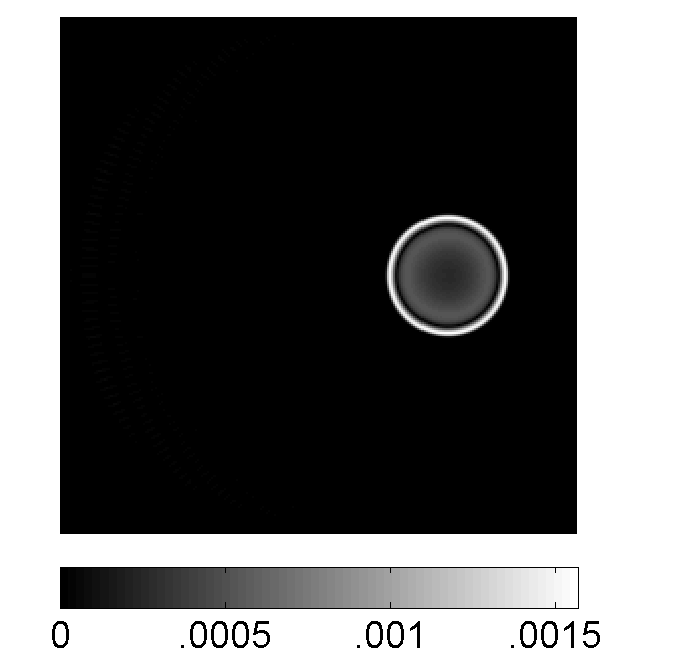}
  \includegraphics[trim = 20 0 20 0, clip, width=0.23\textwidth]{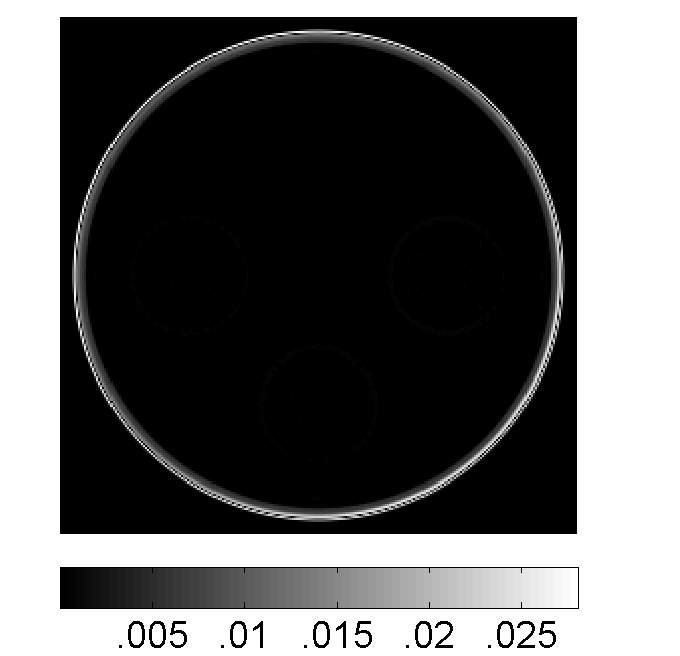} \\
  \includegraphics[trim = 0 0 0 0, clip, width=0.46\textwidth]{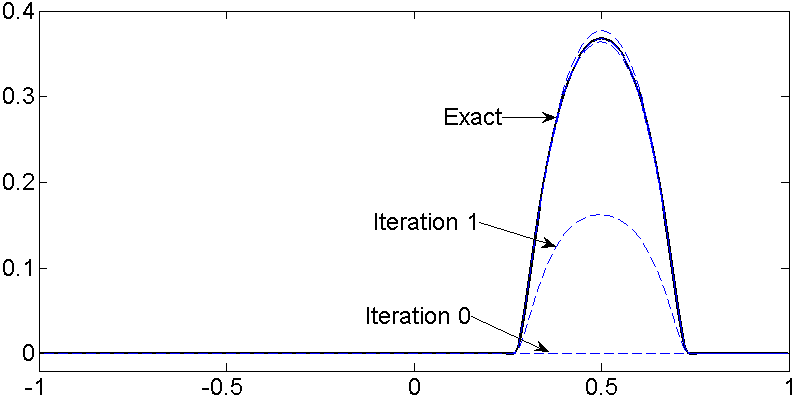}\;\;\;\;
  \includegraphics[trim = 0 0 0 0, clip, width=0.46\textwidth]{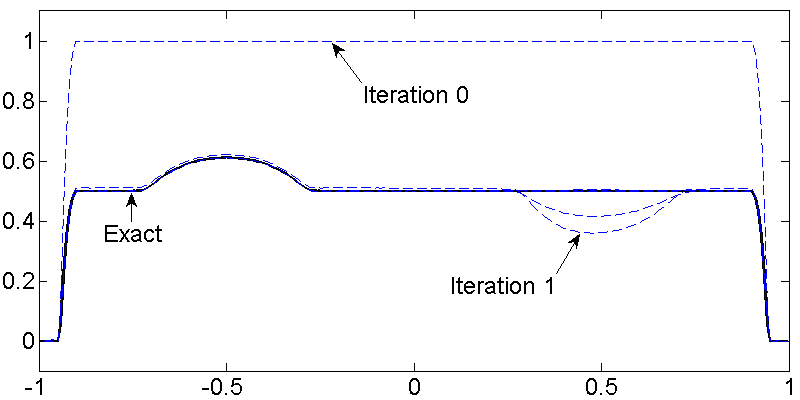}
  \caption{Top row: pointwise error on $(a_1,f_1)$ (left) and $(a_2,f_2)$ (right) after convergence of the Newton-Raphson algorithm. Bottom row: for the reconstruction of $a_2$ (left) and $f_2$ (right), cut plots of the iterations at $\{x=0\}$.}
  \label{fig:Exp12recons}
\end{figure}

Pointwise errors on reconstructions are shown on Fig. \ref{fig:Exp12recons}, where in both scenarios, the reconstruction is excellent. The cut plots in the second row of Fig. \ref{fig:Exp12recons} illustrate the speed of convergence of the method as well as the cross-talk between the reconstructed quantities during the iterations.

\subsection{Reconstruction of non-smooth coefficients and robustness to noise} \label{sec:nonsmooth}
\begin{figure}[htpb]
  \centering
  \includegraphics[height=0.18\textheight]{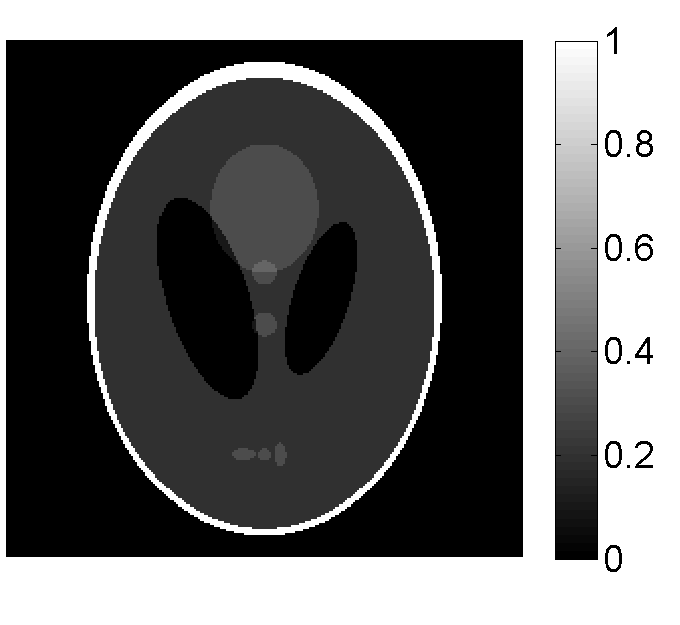}\;\;\;
  \includegraphics[height=0.18\textheight]{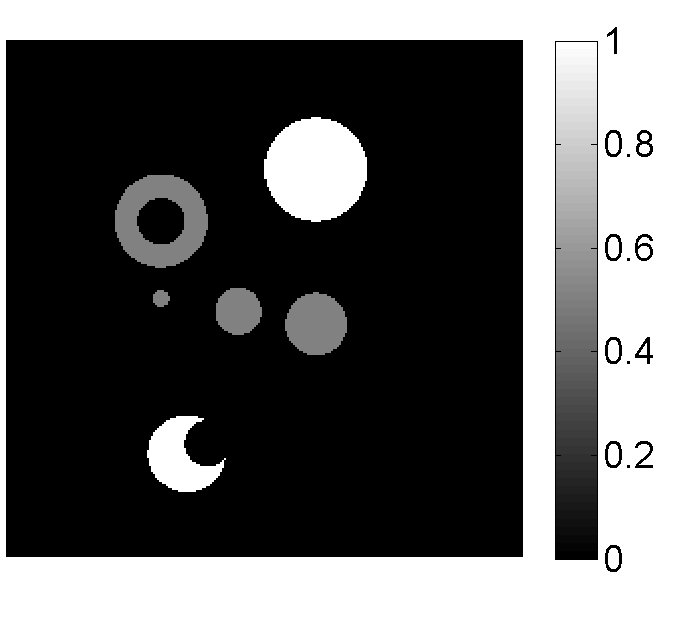}
  \caption{Examples of discontinuous coefficients $a$ (left) and $f$ (right).}
  \label{fig:discont}
\end{figure}

We now consider the case of discontinuous unknown coefficients. We run three simulations using the same discontinuous unknowns (shown in Fig. \ref{fig:discont}), one using noiseless data and the other two polluted with instrumental noise with different levels. 

\paragraph{Noise model.} We add to our measurement a noise of two natures:
\begin{enumerate}
  \item The first kind, modelling instrumental noise, is characterized by an amplitude $A$ so that, each data pixel value $p$ is replaced by a draw ``$A\cdot$Pois$(\frac{p}{A})$''.
  \item After this is done, a background noise is added, characterized by a bias value \[B = \frac{\#\text{added background photons}}{\#\text{photons measured}}. \]
    After deciding a value for a quantum $q$ of energy representing one photon, for each additional photon, we add $q$ to a pixel chosen at random with uniform probability among all data pixels. 
\end{enumerate}

\begin{figure}[htpb]
  \centering
  \includegraphics[trim = 0 25 0 25, clip, height=0.18\textheight]{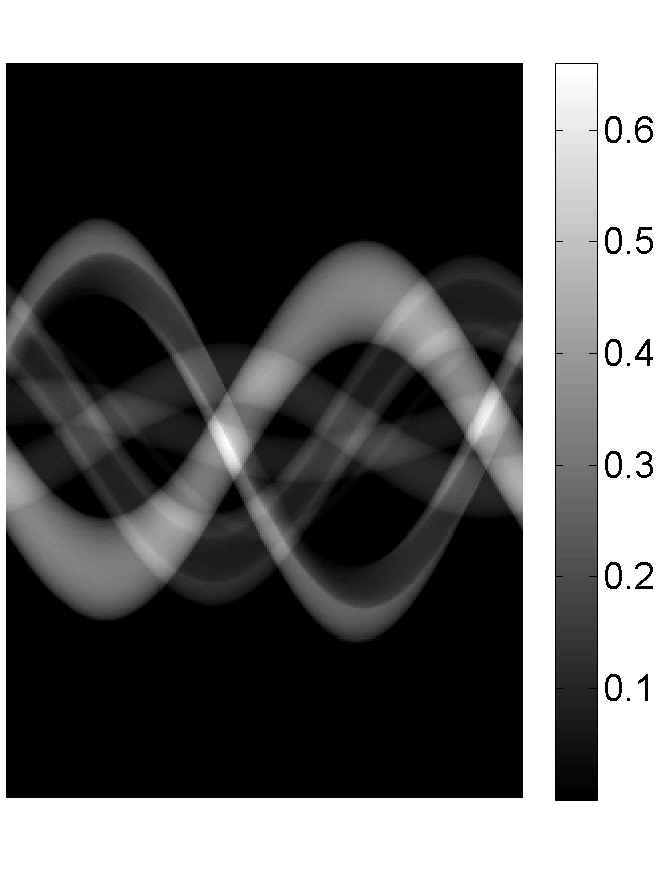}
  \includegraphics[trim = 0 25 0 25, clip, height=0.18\textheight]{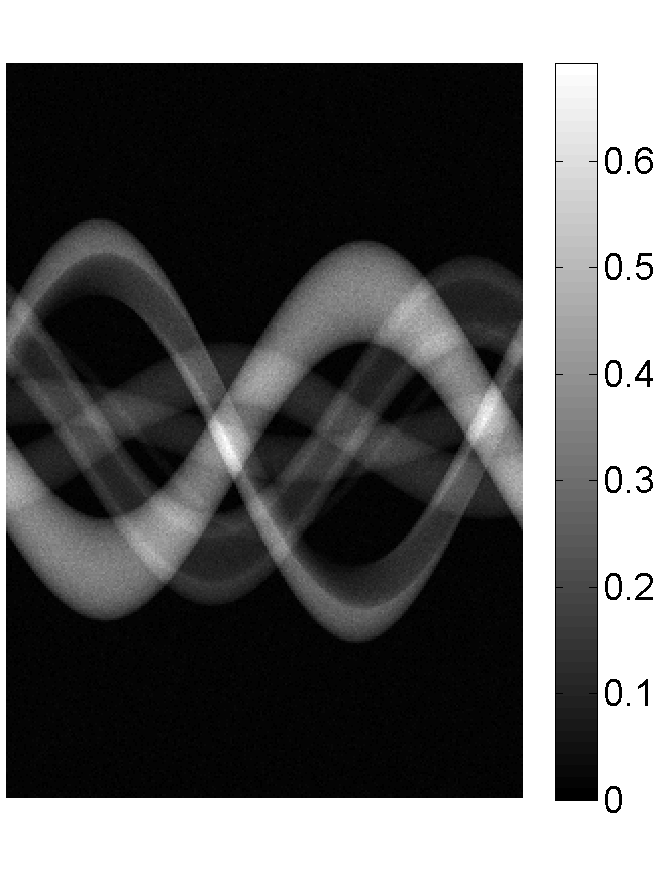}
  \includegraphics[trim = 0 25 0 25, clip, height=0.18\textheight]{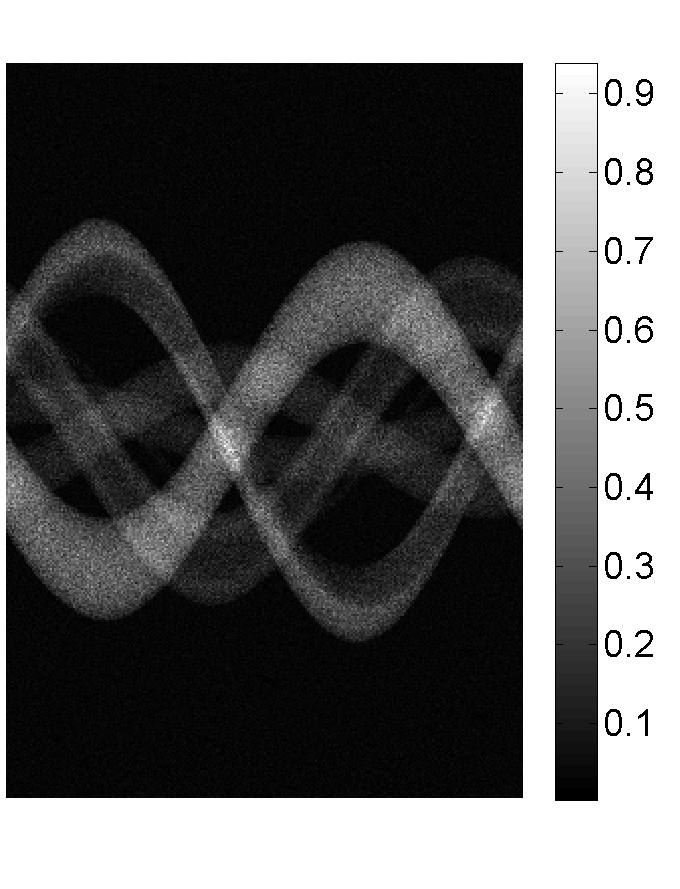}  \\
  \includegraphics[trim = 0 25 0 25, clip, height=0.18\textheight]{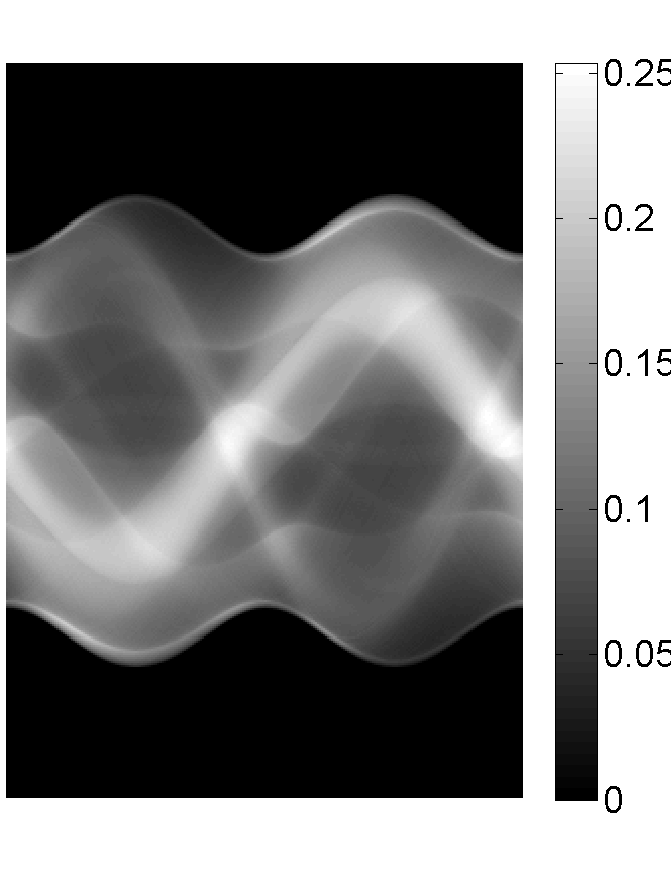}
  \includegraphics[trim = 0 25 0 25, clip, height=0.18\textheight]{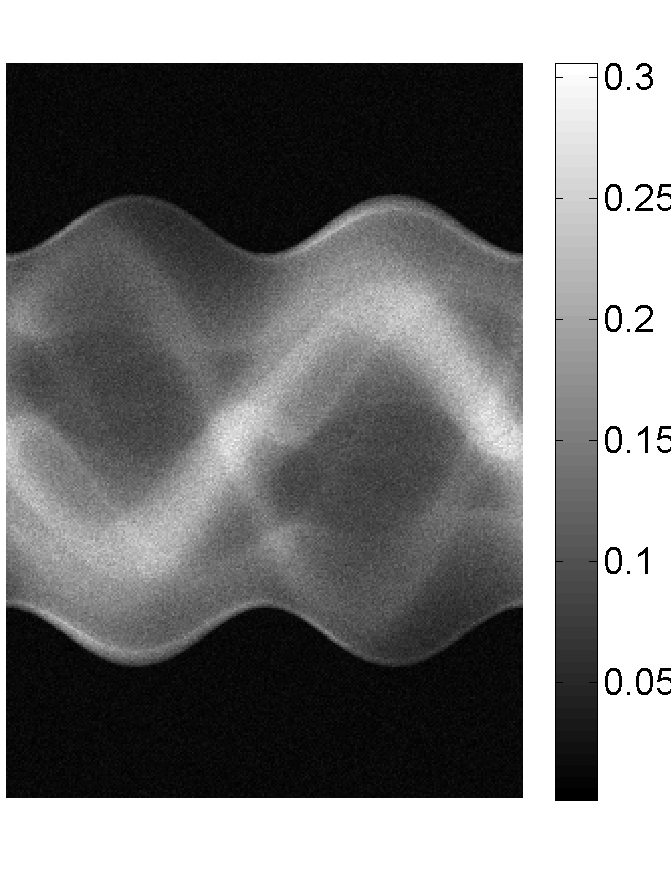}
  \includegraphics[trim = 0 25 0 25, clip, height=0.18\textheight]{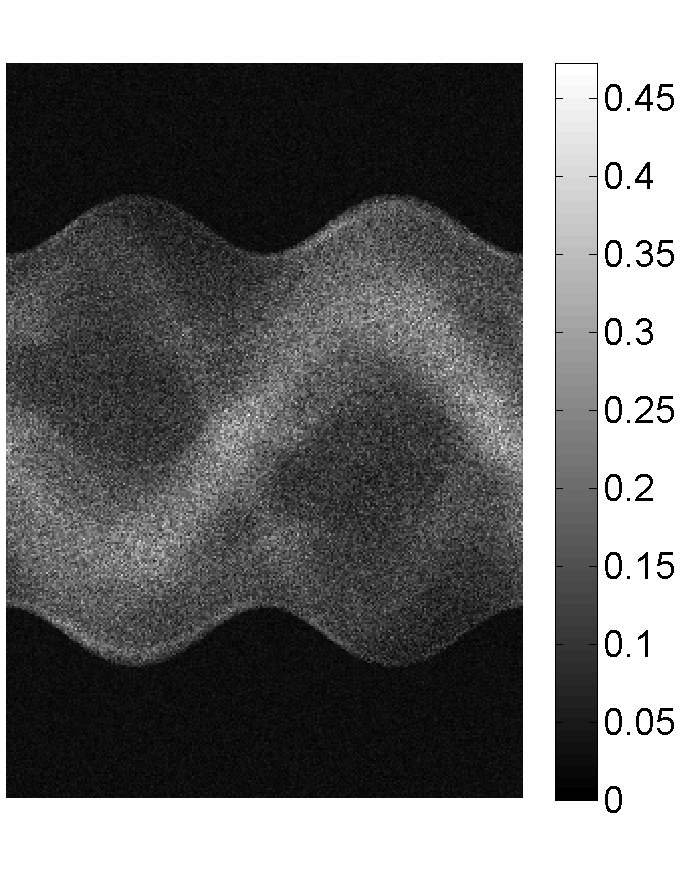}    
  \caption{Forward data $\A_0(a,f)$ (top row) and $\A_1(a,f)$ (bottom row), with $(a,f)$ given in Fig. \ref{fig:discont}. Left to right: noiseless, low noise, high noise.}
  \label{fig:datanoise}
\end{figure}

\begin{figure}[htpb]
  \centering
  \includegraphics[trim = 30 0 30 0, clip, width=0.24\textwidth]{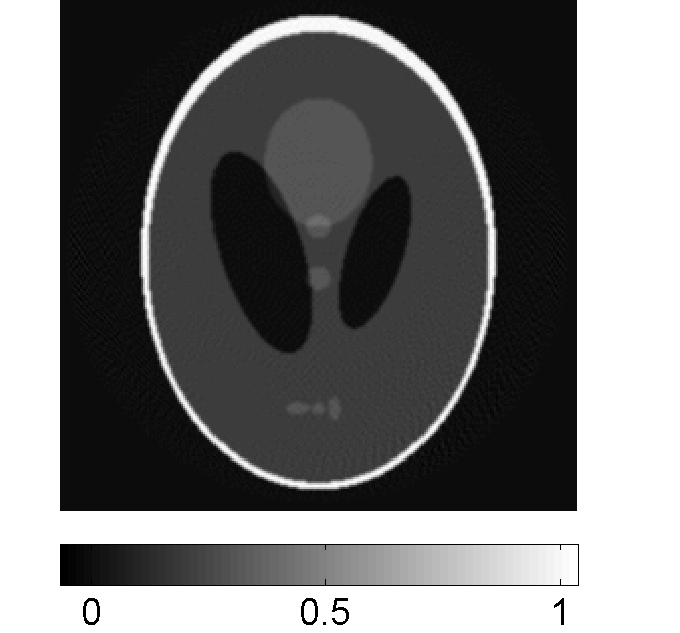}
  \includegraphics[trim = 30 0 30 0, clip, width=0.24\textwidth]{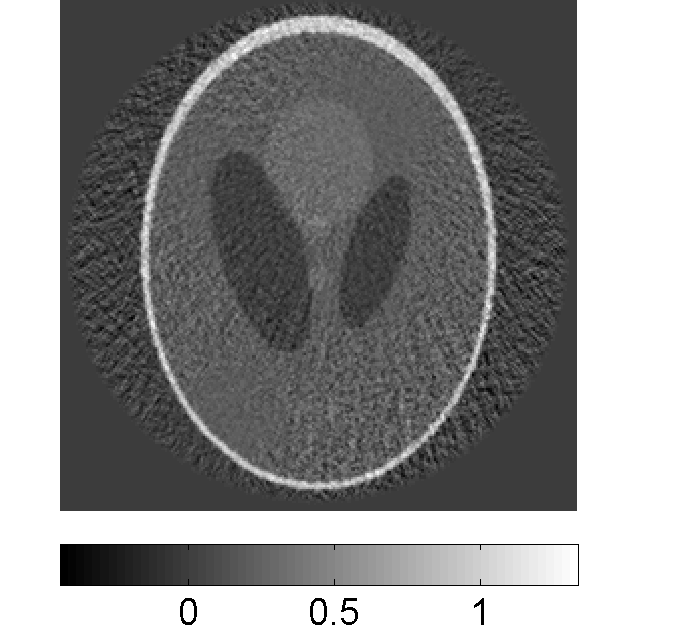}
  \includegraphics[trim = 30 0 30 0, clip, width=0.24\textwidth]{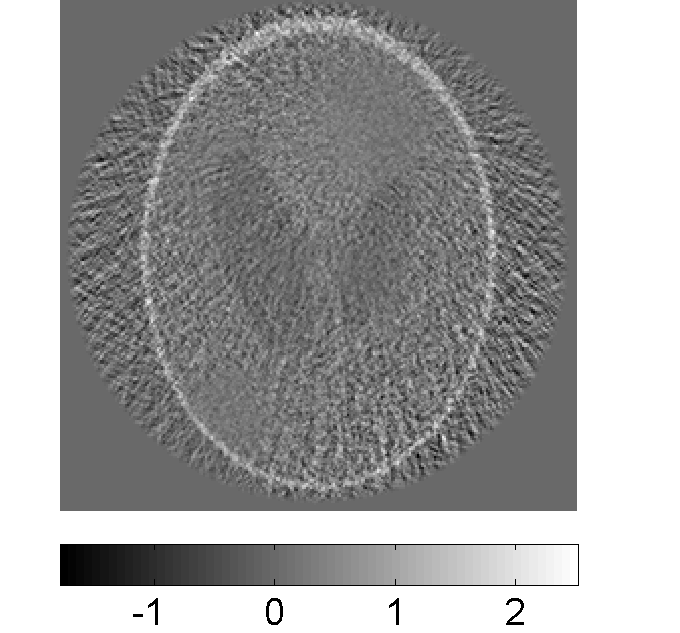}
  \includegraphics[trim = 20 0 30 0, clip, width=0.24\textwidth]{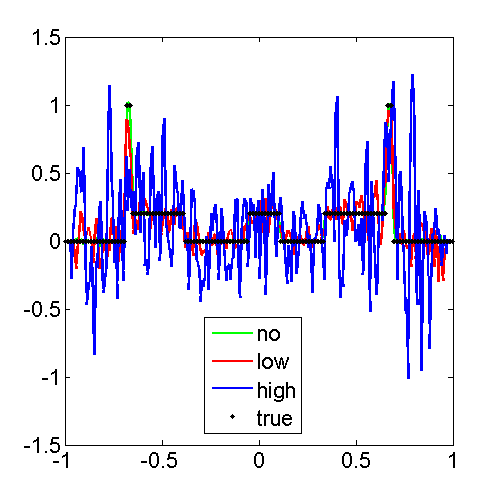}\\
  \includegraphics[trim = 30 0 30 0, clip, width=0.24\textwidth]{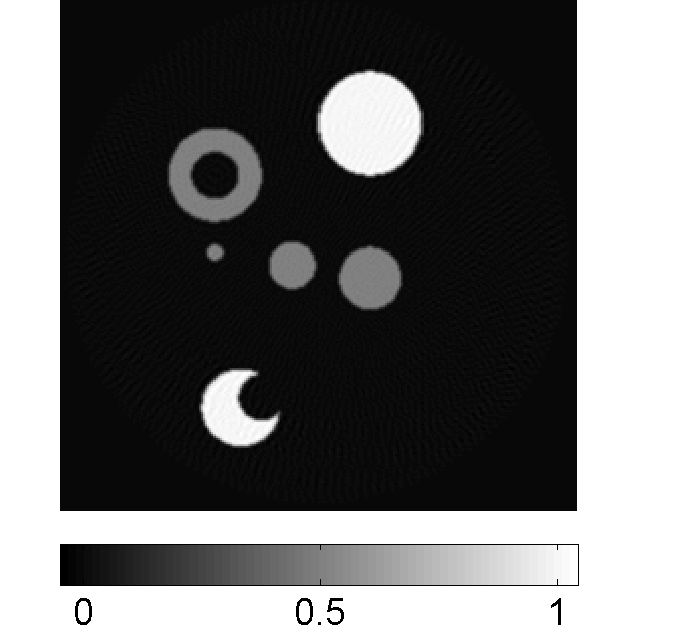}
  \includegraphics[trim = 30 0 30 0, clip, width=0.24\textwidth]{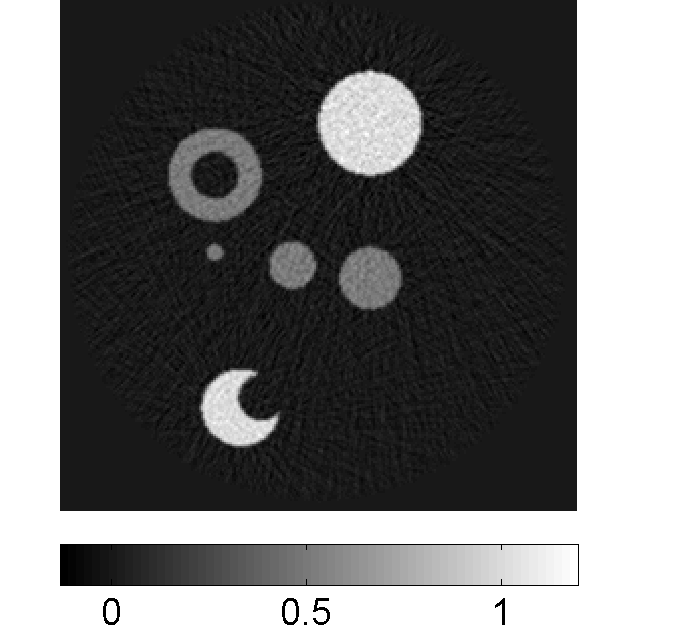}
  \includegraphics[trim = 30 0 30 0, clip, width=0.24\textwidth]{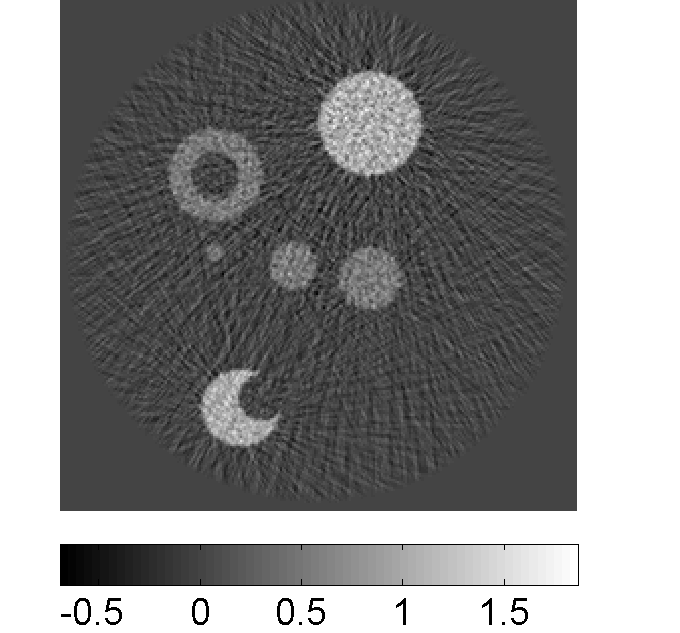}  
  \includegraphics[trim = 20 0 30 0, clip, width=0.24\textwidth]{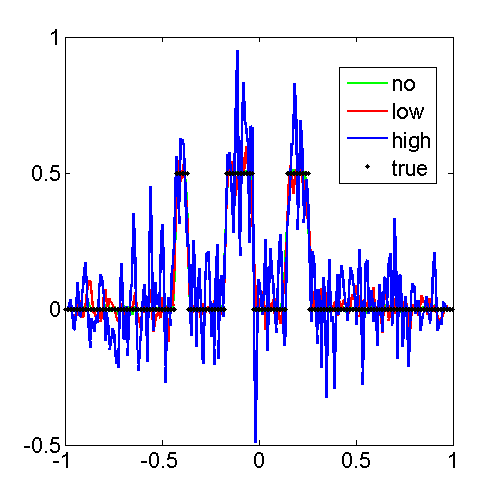}  
  \caption{Reconstructed $a$ (top row) and $f$ (bottom row) after convergence. Left to right: noiseless, low noise, high noise, cut plots at $\{x=0\}$.}
  \label{fig:disconterrors}
\end{figure}

\noindent The experiments with ``low noise'' and ``high noise'' below are carried out with the respective values $(A,B) = (0.2,0.5)$ and $(A,B) = (0.4,5)$. The forward data $(\A_0,\A_1)$ are displayed on Fig. \ref{fig:datanoise}, and the errors after convergence in all three cases (noiseless, low noise, high noise) are displayed in Fig. \ref{fig:disconterrors}. The relative mean square (RMS) errors after 8 iterations are summarized in Table \ref{tab:1}.

\begin{table}
    \centering
    \begin{tabular}[c]{l||r|r|r}
	& noiseless & low noise & high noise \\
	\hline\hline
	RMS on $a$ & 0.2\% & 38.6\% & 127.3\% \\
	\hline
	RMS on $f$ & 0.13\% & 18.7\% & 55.1\%
    \end{tabular}
    \caption{Relative Mean-Square errors on $a$ and $f$ at 8 iterations corresponding to the plots displayed on Fig. \ref{fig:disconterrors}. }
    \label{tab:1}
\end{table}

\paragraph{Comments.} 
\begin{enumerate}
    \item Regarding the choice of initial guess, $f$ should be chosen at first to be a non-vanishing function, so as to prevent the vanishing of the focused transform $M[a,f]$ which appears in denominators of subsequent operations. As seen above, the choice $f\equiv 1$ leads to satisfactory convergence.
    \item As may be seen on Fig. \ref{fig:disconterrors}, although strong additive noise impacts the reconstructions badly, one may notice on the cut plots that the oscillations on reconstructions average about each constant value, leading us to believe that a penalization term (e.g. total variation norm) favoring piecewise constant functions  would re-establish good convergence. Additionally, noise in data may make the algorithm give negative values to both $a$ and $f$, although both quantities are physically nonnegative. This may be avoided by introducing at each iteration a projection step onto nonnegative functions (i.e. of the form $a(x)=\max(a(x),0)$), at the cost of losing the first property of averaging around the correct constant values. Finding an algorithm taking additive noise into account while respecting physics-based criteria appropriately will be the object of future work. 
\end{enumerate}


\section*{Acknowledgments} M.C. was partially funded by Conicyt-Chile grant Fondecyt \#1141189. FM was partially funded by NSF grant No. 1265958. A.O. was partially funded by Conicyt-Chile grants Fondecyt \#1110290 and Conicyt ACT1106.

\section*{References}
\bibliographystyle{plain}

\bibliography{bibliografia}

\end{document}